\documentclass[11pt]{article}
\usepackage{amsmath}
\usepackage{amssymb}
\usepackage{amsfonts}
\usepackage{tikz}
\usepackage{graphicx}
\usepackage{float}
\usepackage[affil-it]{authblk}

\newtheorem{theorem}{Theorem}

\newtheorem{lemma}[theorem]{Lemma}

\newcommand{\func}[1]{\text{#1}}

\title{A Mixed Variational Formulation for the Wellposedness and Numerical Approximation of a PDE Model Arising in a 3-D Fluid-Structure Interaction}

\author{George Avalos
\thanks{e-mail: \texttt{gavalos@math.unl.edu}}}
\author{Thomas J. Clark
\thanks{e-mail: \texttt{s-tclark15@math.unl.edu}}}
\affil{Department of Mathematics,\\ University of Nebraska-Lincoln}

\begin{document}

\maketitle

\begin{abstract}
We will present qualitative and numerical results on a partial differential equation (PDE) system which models a certain fluid-structure
dynamics.  The wellposedness of this PDE model is established by means of constructing for it a nonstandard semigroup generator representation;
this representation is essentially accomplished by an appropriate elimination of the pressure. This coupled PDE model involves the Stokes system
which evolves on a three dimensional domain $\mathcal{O}$ being coupled to a fourth order plate equation, possibly with rotational inertia
parameter $\rho >0$, which evolves on a flat portion $\Omega$ of the boundary of $\mathcal{O}$. The coupling on $\Omega$ is implemented via the
Dirichlet trace of the Stokes system fluid variable - and so the no-slip condition is necessarily not in play - and via the Dirichlet boundary
trace of the pressure, which essentially acts as a forcing term on this elastic portion of the boundary. We note here that inasmuch as the
Stokes fluid velocity does not vanish on $\Omega$, the pressure variable cannot be eliminated by the classic Leray projector; instead, the
pressure is identified as the solution of a certain elliptic boundary value problem.  Eventually, wellposedness of this fluid-structure dynamics
is attained through a certain nonstandard variational (``inf-sup") formulation.  Subsequently we show how our constructive proof of
wellposedness naturally gives rise to a certain mixed finite element method for numerically approximating solutions of this fluid-structure
dynamics.

\end{abstract}

\medskip

\section{The PDE and Setting for Wellposedness}

One of our main objectives in this work is to provide a proof for semigroup wellposedness with respect to the fluid-structure partial
differential equation (PDE) model considered in \cite{igor} - see also \cite{Chambolle} and \cite{igor11}. The proof here will be wholly
different than that originally given in \cite{igor}, and has the virtue of giving insight into a mixed finite element method (FEM) formulation
so as to numerically approximate the solution of the fluid and structure variables.  A numerical analysis involving this fluid-structure FEM
will constitute the second part of this work. Throughout, we will consider situations in which either the \textquotedblleft
Euler-Bernoulli\textquotedblright\ or \textquotedblleft Kirchhoff\textquotedblright\ plate PDE is in place to describe the structural component
of the fluid-structure model (only the Euler-Bernoulli is considered in \cite{igor}). The geometrical situation will be identical to
that in \cite{igor}. We state it here verbatim: $\mathcal{O}\subset \mathbb{R%
}^{3}$ will be a bounded domain with sufficiently smooth boundary. Moreover,
$\partial \mathcal{O}=\bar{\Omega}\cup \bar{S}$, with $\Omega \cap
S=\varnothing $, and specifically%
\begin{equation*}
\Omega \subset \left\{ x=(x_{1,}x_{2},0)\right\} \text{, and surface }%
S\subset \left\{ x=(x_{1,}x_{2},x_{3}):x_{3}\leq 0\right\} .
\end{equation*}%
In consequence, if $\nu (x)$ denotes the exterior unit normal vector to $%
\partial \mathcal{O}$, then
\begin{equation}
\left. \nu \right\vert _{\Omega }=\left[ 0,0,1\right] .  \label{normal}
\end{equation}

\begin{figure}[H]
\begin{center}
\begin{tikzpicture}[scale=1.3]
\draw[left color=black!10,right color=black!20,middle color=black!50, ultra thick] (-2,0,0) to [out=0, in=180] (2,0,0) to [out=270, in = 0]
(0,-3,0) to [out=180, in =270] (-2,0,0);

\draw [fill=black!60, ultra thick] (-2,0,0) to [out=80, in=205](-1.214,.607,0) to [out=25, in=180](0.2,.8,0) to [out=0, in=155] (1.614,.507,0)
to [out=335, in=100](2,0,0) to [out=270, in=25] (1.214,-.507,0) to [out=205, in=0](-0.2,-.8,0) [out=180, in=335] to (-1.614,-.607,0) to
[out=155, in=260] (-2,0,0);

\draw [dashed, thin] (-1.7,-1.7,0) to [out=80, in=225](-.6,-1.3,0) to [out=25, in=180](0.35,-1.1,0) to [out=0, in=155] (1.3,-1.4,0) to [out=335,
in=100](1.65,-1.7,0) to [out=270, in=25] (0.9,-2.0,0) to [out=205, in=0](-0.2,-2.2,0) [out=180, in=335] to (-1.514,-2.0) to [out=155, in=290]
(-1.65,-1.7,0);

\node at (0.2,0.1,0) {{\LARGE$\Omega$}};

\node at (1.95,-1.5,0) {{\LARGE $S$}};

\node at (-0.3,-1.6,0) {{\LARGE $\mathcal{O}$}};
\end{tikzpicture}
\caption{The Fluid-Structure Geometry}
\end{center}
\end{figure}
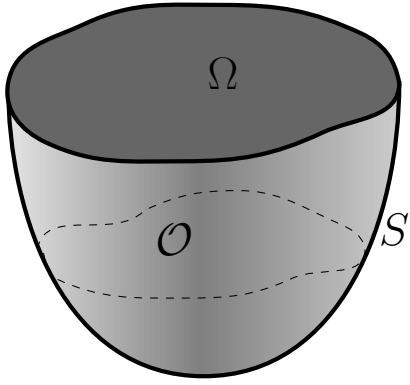

With \textquotedblleft rotational inertia parameter\textquotedblright\ $\rho
\geq 0$, the PDE model is as follows, in solution variables $w(x,t)$, $%
u(x,t)=[u^{1}(x,t),u^{2}(x,t),u^{3}(x,t)]$ and $p(x,t)$:
\begin{eqnarray}
&&w_{tt}-\rho \Delta w_{tt}+\Delta ^{2}w=\left. p\right\vert _{\Omega }\text{
\ in }\Omega \times (0,T);  \label{1} \\
&&w=\frac{\partial w}{\partial \nu }=0\text{ \ on }\partial \Omega ;
\label{2} \\
&&u_{t}-\Delta u+\nabla p=0\text{ \ in }\mathcal{O}\times (0,T);  \label{3}
\\
&&\mathrm{div}(u)=0\text{ \ in }\mathcal{O}\times (0,T);  \label{4} \\
&&u=\vec{0}\text{ on }S\text{ \ and }u=[u^{1},u^{2},u^{3}]=[0,0,w_{t}]\text{ \ on }%
\Omega ,  \label{5}
\end{eqnarray}%
with initial conditions%
\begin{equation}
\lbrack w(0),w_{t}(0),u(0)]=[w_{1},w_{2},u_{0}]\in \mathbf{H}_{\rho }\text{.} \label{ic}
\end{equation}%
(So when $\rho = 0$, Euler-Bernoulli plate dynamics are in play; when $\rho >0$ we have instead the Kirchhoff plate.)  Here, the space of initial data $\mathbf{H}_{\rho }$ is defined as follows:
Let
\begin{equation}
\mathcal{H}_{\mathrm{fluid}}=\left\{ f\in \mathbf{L}^{2}(\mathcal{O}):%
\mathrm{div}(f)=0\text{; }\left. f\cdot \nu \right\vert _{S}=0\right\} ;
\label{H_f}
\end{equation}%
and%
\begin{equation}
W_{\rho }=\left\{
\begin{array}{l}
\dfrac{L^{2}(\Omega )}{\mathbb{R}}\text{, \ if }\rho =0, \\
\\
\dfrac{H_{0}^{1}(\Omega )}{\mathbb{R}}\text{, \ if }\rho
>0.%
\end{array}%
\right.  \label{velocity}
\end{equation}

Therewith, we then set
\begin{eqnarray}
\mathbf{H}_{\rho } &=&\Big\{\left[ \omega_{1},\omega_{2},f\right] \in \left[ H_{0}^{2}(\Omega )\cap \frac{L^{2}(\Omega )}{\mathbb{R}}\right]
\times
W_{\rho }\times \mathcal{H}_{\mathrm{fluid}}  \notag \\
&&\text{ \ \ \ \ }\text{with }\left. f\cdot \nu \right\vert _{\Omega }=[0,0,f^{3}]\cdot \lbrack 0,0,1]=\omega_{2}\Big\}.  \label{energy}
\end{eqnarray}

Moreover, let $A_{D}:L^{2}(\Omega )\rightarrow L^{2}(\Omega )$ be
given by%
\begin{equation}
A_{D}g=-\Delta g\text{, \ \ }D(A_{D})=H^{2}(\Omega )\cap H_{0}^{1}(\Omega ).
\label{dirichlet}
\end{equation}%
If we subsequently make the denotation for all $\rho \geq 0$,%
\begin{equation}
P_{\rho }=I+\rho A_{D}\text{, \ }D(P_{\rho })=\left\{
\begin{array}{l}
L^{2}(\Omega )\text{, \ if }\rho =0, \\
D(A_{D})\text{, \ if }\rho >0,%
\end{array}%
\right.  \label{P}
\end{equation}%
then the mechanical PDE component (\ref{1})-(\ref{2}) can be written as
\begin{equation*}
P_{\rho }w_{tt}+\Delta ^{2}w=\left. p\right\vert _{\Omega }\text{ on\ }(0,T).
\end{equation*}%
Using the characterization from \cite{grisvard} that%
\begin{equation*}
\text{\ }D(P_{\rho }^{\frac{1}{2}})=\left\{
\begin{array}{l}
L^{2}(\Omega )\text{, \ if }\rho =0, \\
H_{0}^{1}(\Omega )\text{, \ if }\rho >0,%
\end{array}%
\right.
\end{equation*}%
then we can endow the Hilbert space $\mathbf{H}_{\rho }$ with norm-inducing
inner product%
\begin{equation*}
\left( \big[ \omega_{1},\omega_{2},f\big] ,\big[ \tilde{\omega}_{1},\tilde{\omega}%
_{2},\tilde{f}\big] \right) _{\mathbf{H}_{\rho }}=(\Delta \omega_{1},\Delta
\tilde{\omega}_{1})_{\Omega }+(P_{\rho }^{\frac{1}{2}}\omega_{2},P_{\rho }^{\frac{%
1}{2}}\tilde{\omega}_{2})_{\Omega }+(f,\tilde{f})_{\mathcal{O}},
\end{equation*}%
where $(\cdot ,\cdot )_{\Omega }$ and $(\cdot ,\cdot )_{\mathcal{O}}$ are
the $L^{2}$-inner products on their respective geometries.

We note here, as there was in \cite{igor}, the necessity for imposing that wave
initial displacement and velocity each have zero mean average. To see this:
Invoking the boundary condition (\ref{5}) and the fact that normal vector $%
\nu =[0,0,1]$ on $\Omega $, we have then by Green's formula, that for all $%
t\geq 0$,%
\begin{equation}
\int_{\Omega }w_{t}(t)d\Omega =\int_{\Omega }u^{3}(t)d\Omega =\int_{\partial
\mathcal{O}}u(t)\cdot \nu d\partial \mathcal{O}=0.  \label{zero}
\end{equation}%
And so we have necessarily,%
\begin{equation*}
\int_{\Omega }w(t)d\Omega =\int_{\Omega }w_{1}d\Omega \text{, for all }t\geq 0\text{.}
\end{equation*}%
This accounts for the choice of the structural finite energy space
components for $\mathbf{H}_{\rho }$, in (\ref{energy}).

As we said, our proof of wellposedness hinges upon demonstrating the existence of a modeling $C_0$-semigroup $\left\{ e^{\mathcal{A}_{\rho
}t}\right\} _{t\geq 0}\subset \mathcal{L}(\mathbf{H}_{\rho })$, for appropriate generator $\mathcal{A}_\rho :\mathbf{H}_{\rho }\rightarrow
\mathbf{H}_{\rho }$. Subsequently, by means of this family, the solution to (\ref{1})-(\ref{ic}), for initial data $[w_{1},w_{2},u_{0}]\in
\mathbf{H}_{\rho }$, will then of course be given via the relation
\begin{equation}
\left[
\begin{array}{c}
w(t) \\
w_{t}(t) \\
u(t)%
\end{array}%
\right] =e^{\mathcal{A}_{\rho }t}\left[
\begin{array}{c}
w_{1} \\
w_{2} \\
u_{0}%
\end{array}%
\right] \in C([0,T];\mathbf{H}_{\rho })\text{.}  \label{semi}
\end{equation}%
Our particular choice here of generator $\mathcal{A}_{\rho }:\mathbf{H}%
_{\rho }\rightarrow \mathbf{H}_{\rho }$ is dictated by the following
consideration:

If $p(t)$ is a viable pressure variable for (\ref{1})-(\ref{ic}), then
pointwise in time $p(t)$ necessarily satisfies the following boundary value
problem:%
\begin{eqnarray}
&&\Delta p=0\text{ \ in }\mathcal{O}\text{;}  \label{bvp1} \\
&&\frac{\partial p}{\partial \nu }+P_{\rho }^{-1}p=P_{\rho }^{-1}\Delta ^{2}w+\left. \Delta u^{3}\right\vert _{\Omega }\text{ \ on }\Omega;
\label{bvp2} \\
&&\frac{\partial p}{\partial \nu }=\left. \Delta u\cdot \nu \right\vert _{S}%
\text{ \ on }S.  \label{bvp3}
\end{eqnarray}%
To show the validity of (\ref{bvp1})-(\ref{bvp3}): taking the divergence of
both sides of (\ref{3}) and using the divergence free condition in (\ref{4})
yields equation (\ref{bvp1}). Moreover, dotting both sides of (\ref{3}) with
the unit normal vector $\nu $, and then subsequently taking the resulting
trace on $S$ will yield the boundary condition (\ref{bvp3}). (Implicitly, we
are also using the fact that $u=0$ on $S$.)

Finally, we consider the particular geometry which is in play (where $\nu =[0,0,1]$ on $\Omega $ to establish (\ref{bvp2})). Using the equation (\ref{1}) and boundary
condition (\ref{5}), we have on $\Omega $
\begin{align*}
P_{\rho }^{-1}\Delta ^{2}w&=-w_{tt}+\left. P_{\rho }^{-1}p\right\vert
_{\Omega } \\
&=-\frac{d}{dt}(0,0,w_{t})\cdot \nu +\left. P_{\rho }^{-1}p\right\vert
_{\Omega } \\
&=-\left[ u_{t}\cdot \nu \right] _{\Omega }+\left. P_{\rho
}^{-1}p\right\vert _{\Omega } \\
&=-\left[ \Delta u\cdot \nu \right] _{\Omega }+\left. \frac{\partial p}{%
\partial \nu }\right\vert _{\Omega }+\left. P_{\rho }^{-1}p\right\vert
_{\Omega },
\end{align*}%
which gives (\ref{bvp2}).

The BVP (\ref{bvp1})-(\ref{bvp3}) can be solved through the agency of the
following \textquotedblleft Robin\textquotedblright\ maps $R_{\rho}$ and $%
\tilde{R}_{\rho}$: We define%
\begin{eqnarray}
R_{\rho }g &=&f\Leftrightarrow \left\{ \Delta f=0\text{ \ in }\mathcal{O};%
\text{ }\frac{\partial f}{\partial \nu }+P_{\rho }^{-1}f=g\text{ \ on }%
\Omega \text{; \ }\frac{\partial f}{\partial \nu }=0\text{ on }S\right\} .
\label{R1} \\
&&  \notag \\
\tilde{R}_{\rho }g &=&f\Leftrightarrow \left\{ \Delta f=0\text{ \ in }%
\mathcal{O};\text{ }\frac{\partial f}{\partial \nu }+P_{\rho }^{-1}f=0\text{
\ on }\Omega \text{; \ }\frac{\partial f}{\partial \nu }=g\text{ on }%
S\right\} .  \label{R2}
\end{eqnarray}%
Therewith, we have that for all real $s$,%
\begin{equation}
R_{\rho }\in \mathcal{L}\big(H^{s}(\Omega ),H^{s+\frac{3}{2}}(\mathcal{O})%
\big)\text{; \ }\tilde{R}_{\rho }\in \mathcal{L}\big(H^{s}(S),H^{s+\frac{3}{2%
}}(\mathcal{O})\big).  \label{Rs}
\end{equation}%
(See e.g. \cite{L-M}. We are also using implicity the fact that $P_{\rho
}^{-1}$ is positive definite, self-adjoint on $\Omega $, and moreover manifests elliptic regularity.)

Therewith, the pressure variable $p(t)$, as necessarily the solution of (\ref%
{bvp1})-(\ref{bvp3}), can be written pointwise in time as
\begin{equation}
p(t)=G_{\rho ,1}(w(t))+G_{\rho ,2}(u(t)),  \label{p}
\end{equation}%
where
\begin{eqnarray}
G_{\rho ,1}(w) &=&R_{\rho }(P_{\rho }^{-1}\Delta ^{2}w);  \label{G1} \\
G_{\rho ,2}(u) &=&R_{\rho }(\left. \Delta u^{3}\right\vert _{\Omega })+%
\tilde{R}_{\rho }(\left. \Delta u\cdot \nu \right\vert _{S}).  \label{G2}
\end{eqnarray}

These relations suggest the following choice for the generator $\mathcal{A}%
_{\rho }:\mathbf{H}_{\rho }\rightarrow \mathbf{H}_{\rho }$. We set%
\begin{eqnarray}
&&\mathcal{A}_{\rho }\equiv
\begin{bmatrix}
0 & I & 0 \\
-P_{\rho }^{-1}\Delta ^{2}+P_{\rho }^{-1}G_{\rho ,1}\big|_{\Omega } & 0 &
P_{\rho }^{-1}G_{\rho ,2}\big|_{\Omega } \\
-\nabla G_{\rho ,1} & 0 & \Delta -\nabla G_{\rho ,2}%
\end{bmatrix}%
;  \label{A} \\
&&  \notag \\
&&\text{with }D(\mathcal{A}_{\rho })=\big\{\left[ w_{1},w_{2},u\right] \in
\mathbf{\ H}_{\rho }\text{ satisfying}:  \notag \\
&&\text{ \quad \quad (a) }w_{1}\in \mathcal{S}_{\rho }\equiv \left\{
\begin{array}{l}
H^{4}(\Omega )\cap H_{0}^{2}(\Omega )\text{, \ if }\rho =0; \\
H^{3}(\Omega )\cap H_{0}^{2}(\Omega )\text{, \ if }\rho >0;%
\end{array}%
\right. \label{S} \\
&&\text{ \quad \quad (b) }w_{2}\in H_{0}^{2}(\Omega )\text{, }u\in \mathbf{H}^{2}(%
\mathcal{O});  \notag \\
&&\text{ \quad \quad (c) }u=\vec{0}\text{ on }S\text{ and }u=[0,0,w_{2}]\text{ on }\Omega \big\}, \label{domain}
\end{eqnarray}
(c.f. the generator and domain described in an earlier version of \cite{igor} in arXiv:1109.4324.)   (Note that as $\Delta u \in
L^2(\mathcal{O})$ and $\text{div}(\Delta u) = 0$, then by Theorem 1.2, p. 9 in \cite{temam}, we have the trace regularity
\begin{equation}
\Delta u\cdot \nu \big|_{\partial \mathcal{O}} \in H^{-\frac{1}{2}}(\partial \mathcal{O} ); \label{CritReg}
\end{equation}
and so the pressure term
\begin{equation}
p \equiv G_{\rho ,1}(w_{1})+G_{\rho ,2}(u)\in H^{1}(\mathcal{O}). \label{pInH1}
\end{equation}
In short the domain of $\mathcal{A}_\rho: \mathbf{H}_\rho \to \mathbf{H}_\rho$ is well-defined.)

In what follows, we will have need of solution maps for certain
inhomogeneous Stokes flows. To wit: For given $\phi \in H^{\frac{3}{2}%
}(\Omega )$, let $[\tilde{f}(\phi ),\tilde{\pi}(\phi )]\in \mathbf{H}^{2}(%
\mathcal{O})\times H^{1}(\mathcal{O})/\mathbb{R}$ solve
\begin{equation}
\begin{cases}
\lambda \tilde{f}-\Delta \tilde{f}+\nabla \tilde{\pi}=0\text{ \ in }\mathcal{%
O}, &  \\
\mathrm{div}(\tilde{f})=\dfrac{1\cdot \int_{\Omega }\phi \,d\Omega }{\text{meas}(%
\mathcal{O})}\text{ \ in }\mathcal{O}, &  \\
\tilde{f}\big|_{S}=[0,0,0]\text{ \ on }S, &  \\
\tilde{f}\big|_{\Omega }=[0,0,\phi ]\text{ \ on }\Omega . &
\end{cases}
\label{sol1}
\end{equation}

We note that the classic compatibility condition for solvability is satisfied, and that pressure variable $\tilde{\pi}$ is uniquely defined up
to a constant; see e.g., Theorem 2.4, p. 31 of \cite{temam}. Then by Agmon-Douglis-Nirenberg, we have $\left[ \tilde{f},\tilde{\pi}\right] \in
\mathcal{L}(H^{\frac{3}{2}}(\Omega ),\mathbf{H}^{2}(\mathcal{O})\times H^{1}(%
\mathcal{O})/\mathbb{R)}$, with
\begin{equation}
\Vert \tilde{f}(\phi )\Vert _{\mathbf{H}^{2}(\mathcal{O})}+\Vert \tilde{\pi}%
(\phi )\Vert _{\frac{H^{1}(\mathcal{O})}{\mathbb{R}}}\leq C\Vert \phi \Vert _{H^{\frac{3}{2}}(\Omega )}.  \label{map1}
\end{equation}%
(see Proposition 2.2, p. 33 of \cite{temam}).

In a similar way, we define for fluid data $u^{\ast }\in \mathbf{L}^{2}(%
\mathcal{O})$, the solution variables $[\tilde{\mu}(u^{\ast }),\tilde{q}%
(u^{\ast })]\in \mathbf{H}^{2}(\mathcal{O})\times H^{1}(\mathcal{O})/\mathbb{%
R}$, where $[\tilde{\mu},\tilde{q}]$ solve
\begin{equation}
\begin{cases}
\lambda \tilde{\mu}-\Delta \tilde{\mu}+\nabla \tilde{q}=u^{\ast } & \text{in
}\mathcal{O}, \\
\mathrm{div}(\tilde{\mu})=0 & \text{in }\mathcal{O}, \\
\tilde{\mu}\big|_{\partial \mathcal{O}}=\vec{0} & \text{on }\partial
\mathcal{O}.%
\end{cases}
\label{sol2}
\end{equation}%
Again by Agmon-Douglis-Nirenberg we have $\left[ \tilde{\mu},\tilde{q}\right]
\in \mathcal{L}(\mathbf{L}^{2}(\mathcal{O}),\mathbf{H}^{2}(\mathcal{O}%
)\times H^{1}(\mathcal{O})/\mathbb{R)}$, with%
\begin{equation}
\Vert \tilde{\mu}(u^{\ast })\Vert _{\mathbf{H}^{2}(\mathcal{O})}+\Vert \tilde{q}(u^{\ast })\Vert _{\frac{H^{1}(\mathcal{O})}{\mathbb{R}}}\leq
C\Vert u^{\ast }\Vert _{\mathbf{L}^{2}(\mathcal{O})}.  \label{map2}
\end{equation}

These two fluid maps will be invoked in the proof of Theorem \ref{well} below, which will yield solvability of the following resolvent equation
for $\lambda >0$:  Namely, for given $[w_1^*, w_2^*, u^*] \in \mathbf{H}_\rho$, $[w_1, w_2, u] \in D(\mathcal{A}_\rho)$ solves
\begin{equation} \label{resolvent}
(\lambda I - \mathcal{A}_\rho) \begin{bmatrix} w_1\\w_2\\u \end{bmatrix} = \begin{bmatrix} w_1^* \\ w_2^* \\ u^* \end{bmatrix}.
\end{equation}
After setting the pressure variable $p=G_{\rho, 1}(w_{1})+G_{\rho, 2}(u)$ in (\ref{A}), the resolvent equation (\ref{resolvent}) is equivalent to the
following system:
\begin{alignat}{3}
\lambda w_{1}-w_{2}& =w_{1}^{\ast } &\quad& && \text{in }\Omega  ;  \label{w2} \\
\lambda w_{2}+P_{\rho }^{-1}\Delta ^{2}w_{1}-P_{\rho }^{-1}p\big|_{\Omega }&
=w_{2}^{\ast }   &\quad& && \text{in }\Omega  ;  \label{w1} \\
w_{1}=\frac{\partial w_{1}}{\partial \nu }& =0   &\quad& &&  \text{on }\partial \Omega  ; \\
\lambda u-\Delta u+\nabla p& =u^{\ast }   &\quad& && \text{in }\mathcal{O}  ;  \label{ustar}\\
\mathrm{div}(u)& =0  &\quad& &&  \text{in }\mathcal{O} ;  \\
u\big|_{S}& =\left[ 0,0,0\right]  &\quad& &&  \text{on }S  ;  \\
u\big|_{\Omega }& =\left[ 0,0,w_{2}\right]  &\quad& &&  \text{in }\Omega .  \label{end}
\end{alignat}
In particular, it will be seen in Theorem \ref{well} that the solution variable $w_1$ in (\ref{resolvent}) can be recovered through finding the
unique solution pair $[w_1, \tilde{c}] \in H^2_0(\Omega) \times \mathbb{R}$ which solves
\begin{equation}
\begin{cases}
a_{\lambda }(w_{1},\phi )+b(\phi ,\tilde{c})=\mathbb{F}(\phi ) & \forall
\phi \in H_{0}^{2}(\Omega ), \\
b(w_{1},r)=0 & \forall r\in \mathbb{R} ; %
\end{cases}%
  \label{BB1}
\end{equation}%
where:
\begin{align}
a_{\lambda }(\psi ,\phi )& =\lambda ^{2}(P_{\rho }^{1/2}\psi ,P_{\rho }^{1/2}\phi )_{\Omega }+(\Delta \psi ,\Delta \phi )_{\Omega
}+\lambda (\nabla \tilde{f}(\psi ),\nabla \tilde{f}(\phi ))_{\mathcal{O}%
}+\lambda ^{2}(\tilde{f}(\psi ),\tilde{f}(\phi ))_{\mathcal{O}}, \notag \\
& \hspace{0.15in}\forall \,\psi \text{ and }\phi \in H_{0}^{2}(\Omega ); \notag
\\
b(\phi ,r)& =-r\int_{\Omega }\phi \,d\Omega ,\hspace{0.2in}\forall
\phi \in H_{0}^{2}(\Omega )\text{ and }r\in \mathbb{R}; \label{abF}\\
\mathbb{F}(\phi )& =(\nabla \tilde{f}(w_{1}^{\ast }),\nabla \tilde{%
f}(\phi ))_{\mathcal{O}}+\lambda (\tilde{f}(w_{1}^{\ast }),\tilde{f}(\phi
))_{\mathcal{O}}-(\nabla \tilde{\mu}(u^{\ast }),\nabla \tilde{f}(\phi ))_{%
\mathcal{O}} \notag \\
& \hspace{.2in}-\lambda (\tilde{\mu}(u^{\ast }),\tilde{f}(\phi ))_{\mathcal{O}%
}+(u^{\ast },\tilde{f}(\phi ))_{\mathcal{O}}+\left( P_{\rho }(\lambda w_{1}^{\ast }+w_{2}^{\ast }),\phi \right)_{\Omega }, \hspace{0.2in}\forall
\phi \in H_{0}^{2}(\Omega ). \notag
\end{align}

\begin{theorem} \hspace{.2in}
\label{well}
\begin{enumerate}
\item[(i)] The operator $\mathcal{A}_{\rho }:\mathbf{H}_{\rho }\rightarrow \mathbf{H}_{\rho }$ is maximal dissipative.  Therefore by the
Lumer-Phillips Theorem it generates a $C_{0}$-semigroup of contractions $\left\{ e^{\mathcal{A}_{\rho}t}\right\} _{t\geq 0}$ on
$\mathbf{H}_{\rho }$.
\item[(ii)] Let $\lambda >0$ and $[w_1^*, w_2^*, u^*] \in \mathbf{H}_\rho$ be given. (By part (i), there exists
$[w_1, w_2, u] \in D(\mathcal{A}_\rho)$ which solves (\ref{resolvent}).)  Then the structural solution component $w_1$ and constant component
$\tilde{c}$ of the associated pressure term $p$ can be characterized as the solution pair of the variational system (\ref{BB1}). Subsequently
the remaining unknown terms are given by
\begin{align}
w_2 &= \lambda w_1 - w_1^*, \notag \\
u&=\tilde{f}(\lambda w_{1}-w_{1}^{\ast })+\tilde{\mu}(u^{\ast }), \label{up} \\
p&=\tilde{\pi}(\lambda w_{1}-w_{1}^{\ast })+\tilde{q}(u^{\ast })+\tilde{c}, \notag %
\end{align}
(after utilizing the maps $\tilde{f}$ and $\tilde{\pi}$ in (\ref{sol1}) and $\tilde{\mu}$ and $\tilde{q}$ in (\ref{sol2})).
\end{enumerate}
\end{theorem}

\section{Proof of Theorem \protect\ref{well}}

\subsection{Proof of Dissipativity}

Let $\left[ w_{1},w_{2},u\right] \in D(\mathcal{A})$ be given. Therewith, we
set pressure%
\begin{equation}
\pi _{0}\equiv G_{\rho ,1}(w_{1})+G_{\rho ,2}(u).  \label{pi}
\end{equation}%
(So by (\ref{G1})-(\ref{G2}), (\ref{Rs}) and (\ref{pInH1}), we have that $\pi _{0}\in H^{1}(\mathcal{O})$.) We have then, upon using the
definition of the domain in (\ref{domain}),%
\begin{eqnarray}
&&\left( \mathcal{A}_{\rho }\left[
\begin{array}{c}
w_{1} \\
w_{2} \\
u%
\end{array}%
\right] ,\left[
\begin{array}{c}
w_{1} \\
w_{2} \\
u%
\end{array}%
\right] \right) _{\mathbf{H}_{\rho }}  \notag \\
&=&(\Delta w_{2},\Delta w_{1})_{\Omega }+(-\Delta ^{2}w_{1}+\left. \pi
_{0}\right\vert _{\Omega },w_{2})_{\Omega }+\left( \Delta u-\nabla \pi
_{0},u\right) _{\mathcal{O}}  \notag \\
&&  \notag \\
&=&(\Delta w_{2},\Delta w_{1})_{\Omega }+\left( \nabla \Delta w_{1},\nabla
w_{2}\right) _{\Omega }+(\left. \pi _{0}\right\vert _{\Omega
}(0,0,1),(u^{1},u^{2},w_{2}))_{\Omega }  \notag \\
&&\text{ \ }-(\nabla u,\nabla u)_{\mathcal{O}}+\left\langle \frac{\partial u%
}{\partial \nu },u\right\rangle _{\Omega }-\left\langle \pi _{0}\nu
,u\right\rangle _{\Omega }  \notag \\
&&  \notag \\
&=&(\Delta w_{2},\Delta w_{1})_{\Omega }-\left( \Delta w_{1},\Delta
w_{2}\right) _{\Omega }-(\nabla u,\nabla u)_{\mathcal{O}}  \notag \\
&&+\left( \left[
\begin{array}{c}
\partial _{x_{3}}u^{1} \\
\partial _{x_{3}}u^{2} \\
\partial _{x_{3}}u^{3}%
\end{array}%
\right] ,\left[
\begin{array}{c}
0 \\
0 \\
u^{3}%
\end{array}%
\right] \right) _{\Omega }  \notag \\
&=&-2i\func{Im}\left( \Delta w_{1},\Delta w_{2}\right) _{\Omega }-\left\Vert
\nabla u\right\Vert _{\mathcal{O}}^{2},  \label{dissi}
\end{eqnarray}%
where in the last step, we have used we have used $u=\vec{0}$ on $S$, $\mathrm{div}%
(u)=0$ and $[u^{1},u^{2},u^{3}]=[0,0,w_{2}]$ on $\Omega $. This establishes dissipativity.

\subsection{Proof of Maximality}
\label{Maximality}

In what follows we will make use of the Babu\v{s}ka-Brezzi Theorem.  We state it here directly from p. 116 of \cite{kesavan}.
\begin{theorem}[Babu\v{s}ka-Brezzi]
\label{thm:BB} Let $\Sigma$, $M$ be Hilbert spaces and $a:\Sigma \times \Sigma \to \mathbb{R}$, $b:\Sigma \times M \to \mathbb{R}$, bilinear
forms which are continuous.  Let
\begin{equation}
Z=\{\sigma \in \Sigma \, \big| \,  b(\sigma, q) = 0, \,\, \text{ for every } q \in M\}. \label{3119}
\end{equation}
Assume that $a(\cdot, \cdot)$ is $Z$-elliptic, i.e. there exists a constant $\alpha >0$ such that
\begin{equation}
a(\sigma, \sigma) \geq \|\sigma\|^2_\Sigma. \label{3120}
\end{equation}
Assume further that there exists a constant $\beta >0$ such that
\begin{equation}
\sup_{\tau \in \Sigma} \frac{b(\tau, q)}{\|\tau\|_\Sigma} \geq \beta \|q\|_M, \,\, \text{ for every } q \in M. \label{3121}
\end{equation}
Then if $\kappa \in \Sigma$ and $\ell \in M$, there exists a unique pair $(\sigma, p) \in \Sigma \times M$ such that
\begin{align}
a(\sigma, \tau) + b(\tau, p) = (\kappa, \tau), \,\, \text{ for every } \tau \in \Sigma. \label{3122}\\
b(\sigma, q) = (\ell, q), \,\, \text{ for every } q \in M. \label{3123}
\end{align}
\end{theorem}

For $\lambda >0$, we will show that $\text{Range}(\lambda I-\mathcal{A}%
_{\rho })=\mathbf{H}_{\rho }$. To this end, let $[w_{1}^{\ast },w_{2}^{\ast
},u^{\ast }]\in \mathbf{H}_{\rho }$ be given. We must find $%
[w_{1},w_{2},u]\in \mathcal{D}(\mathcal{A}_{\rho })\subset \mathbf{H}_{\rho
} $ which solves
\begin{equation}
(\lambda I-\mathcal{A}_{\rho })%
\begin{bmatrix}
w_{1} \\
w_{2} \\
u%
\end{bmatrix}%
=%
\begin{bmatrix}
w_{1}^{\ast } \\
w_{2}^{\ast } \\
u^{\ast }%
\end{bmatrix}%
.  \label{res}
\end{equation}

Using the structural component (\ref{w2}) and (\ref{w1}), we then have the boundary value problem
\begin{equation}
\begin{cases}
\lambda ^{2}w_{1}+P_{\rho }^{-1}\Delta ^{2}w_{1}-P_{\rho }^{-1}p\big|%
_{\Omega }=\lambda w_{1}^{\ast }+w_{2}^{\ast }\text{ in }\Omega ,\text{ \ } &
\\
w_{1}\big|_{\partial \Omega }=\frac{\partial w_{1}}{\partial \nu }\big|%
_{\partial \Omega }=0\text{ \ on }\partial \Omega , &
\end{cases}
\label{str1}
\end{equation}%
as well as the fluid system
\begin{equation}
\begin{cases}
\lambda u-\Delta u+\nabla p=u^{\ast } & \text{in }\mathcal{O}, \\
\mathrm{div}(u)=0 & \text{in }\mathcal{O}, \\
u\big|_{S}=[0,0,0] & \text{on }S, \\
u\big|_{\Omega }=[0,0,\lambda w_{1}-w_{1}^{\ast }] & \text{on }\Omega .%
\end{cases}
\label{flu1}
\end{equation}

Since $[w_{1},w_{2},u]\in \mathbf{H}_{\rho }$, we also have that
\begin{equation}
\int_{\Omega }w_{1}\,d\Omega =0.  \label{w0}
\end{equation}

With this compatibilty condition in mind, by way of \textquotedblleft
decoupling\textquotedblright\ the systems (\ref{str1}) and (\ref{flu1}), we
proceed as follows. We apply $P_{\rho }$ to the mechanical equation in (\ref%
{str1}) and then multiply by test function $\phi \in H_{0}^{2}(\Omega )$. Subsequently
integrating over $\Omega $ then yields
\begin{equation*}
\left( \lbrack \lambda ^{2}P_{\rho }+\Delta ^{2}]w_{1}-\left. p\right\vert
_{\Omega },\phi \right) _{\Omega }=\left( P_{\rho }(\lambda w_{1}^{\ast
}+w_{2}^{\ast }),\phi \right) _{\Omega }
\end{equation*}%
(where $\left( \cdot ,\cdot \right) _{\Omega }$ might also indicate the
duality pairing between $H_{0}^{2}(\Omega )$ and $H^{-2}(\Omega )$). Upon
integration by parts and using $\displaystyle\phi \big|_{\partial \Omega }=%
\frac{\partial \phi }{\partial \nu }\Big|_{\partial \Omega }=0$, we then
have
\begin{equation}
\lambda ^{2}(P_{\rho }w_{1},\phi )_{\Omega }+(\Delta w_{1},\Delta \phi
)_{\Omega }+\left( \frac{\partial u}{\partial \nu }-p\nu ,%
\begin{bmatrix}
0 \\
0 \\
\phi%
\end{bmatrix}%
\right) _{\Omega }=\left( P_{\rho }(\lambda w_{1}^{\ast }+w_{2}^{\ast
}),\phi \right) _{\Omega }.  \label{wk1}
\end{equation}

(Note that in obtaining this expression we have again used $\mathrm{div}%
(u)=0 $ on $\Omega $ and normal vector $\left. \nu \right\vert _{\Omega
}=[0,0,1]$.)

Now, from (\ref{flu1}) and (\ref{w0}), we can use the fact that in terms of
the maps in (\ref{sol1}) and (\ref{sol2}) above, we can write fluid
variables $u$ and $p$ of (\ref{str1})-(\ref{flu1}) as
\begin{equation}
\begin{cases}
u=\tilde{f}(\lambda w_{1}-w_{1}^{\ast })+\tilde{\mu}(u^{\ast }), \\
p=\tilde{\pi}(\lambda w_{1}-w_{1}^{\ast })+\tilde{q}(u^{\ast })+\tilde{c},%
\end{cases} \label{up56}
\end{equation}%
for some (to be determined) constant $\tilde{c}$, justifying (\ref{up}).

Applying this representation to (\ref{wk1}), subsequently integrating by
parts, and using the mapping in (\ref{sol1}), we then have for every $\phi
\in H_{0}^{2}(\Omega )$,
\begin{align*}
&\lambda ^{2}(P_{\rho }w_{1},\phi )_{\Omega }+(\Delta w_{1},\Delta \phi
)_{\Omega }+\left( \nabla u,\nabla \tilde{f}(\phi )\right) _{\mathcal{O}%
}+\left( \Delta u,\tilde{f}(\phi )\right) _{\mathcal{O}}
-\left( p,\func{div}[\tilde{f}(\phi )]\right) _{\mathcal{O}}-\left( \nabla p,\tilde{f}(\phi )\right) _{\mathcal{O}}\\
& \quad =\left( P_{\rho }(\lambda w_{1}^{\ast }+w_{2}^{\ast }),\phi \right) _{\Omega }.
\end{align*}%
Using now the representations in (\ref{up56}), we have then
\begin{align*}
& \lambda ^{2}(P_{\rho }^{1/2}w_{1},P_{\rho }^{1/2}\phi )_{\Omega }+(\Delta
w_{1},\Delta \phi )_{\Omega }+(\nabla \tilde{f}(\lambda w_{1}-w_{1}^{\ast
}),\nabla \tilde{f}(\phi ))_{\mathcal{O}}+(\nabla \tilde{\mu}(u^{\ast
}),\nabla \tilde{f}(\phi ))_{\mathcal{O}} \\
& +(\Delta \tilde{f}(\lambda w_{1}-w_{1}^{\ast }),\tilde{f}(\phi ))_{%
\mathcal{O}}+(\Delta \tilde{\mu}(u^{\ast }),\tilde{f}(\phi ))_{\mathcal{O}%
}-(\nabla \tilde{\pi}(\lambda w_{1}-w_{1}^{\ast }),\tilde{f}(\phi ))_{%
\mathcal{O}} \\
& -(\nabla \tilde{q}(u^{\ast }),\tilde{f}(\phi ))_{\mathcal{O}}-(\tilde{\pi}%
(\lambda w_{1}-w_{1}^{\ast }),\mathrm{div}(\tilde{f}(\phi )))_{\mathcal{O}}-(%
\tilde{q}(u^{\ast }),\mathrm{div}(\tilde{f}(\phi )))_{\mathcal{O}}-\tilde{c}%
(1,\mathrm{div}(\tilde{f}(\phi )))_{\mathcal{O}} \\
& =\left( P_{\rho }(\lambda w_{1}^{\ast }+w_{2}^{\ast }),\phi \right)
_{\Omega }.
\end{align*}

Now using (\ref{sol1}) and (\ref{sol2}) as well as the fact that $\tilde{\pi}%
,\tilde{q}\in \frac{L^{2}(\mathcal{O})}{\mathbb{R}}$, we rewrite this expression as
\begin{align}
& \lambda ^{2}(P_{\rho }^{1/2}w_{1},P_{\rho }^{1/2}\phi )_{\Omega }+(\Delta
w_{1},\Delta \phi )_{\Omega }+\lambda (\nabla \tilde{f}(w_{1}),\nabla \tilde{%
f}(\phi ))_{\mathcal{O}}+\lambda ^{2}(\tilde{f}(w_{1}),\tilde{f}(\phi ))_{%
\mathcal{O}}-\tilde{c}\int_{\Omega }\phi \,d\Omega  \notag \\
& =(\nabla \tilde{f}(w_{1}^{\ast }),\nabla \tilde{f}(\phi ))_{\mathcal{O}%
}+\lambda (\tilde{f}(w_{1}^{\ast }),\tilde{f}(\phi ))_{\mathcal{O}}-(\nabla
\tilde{\mu}(u^{\ast }),\nabla \tilde{f}(\phi ))_{\mathcal{O}}  \notag \\
& \hspace{0.2in}-\lambda (\tilde{\mu}(u^{\ast }),\tilde{f}(\phi ))_{\mathcal{%
O}}+(u^{\ast },\tilde{f}(\phi ))_{\mathcal{O}}+\left( P_{\rho }(\lambda
w_{1}^{\ast }+w_{2}^{\ast }),\phi \right) _{\Omega }.  \label{var1}
\end{align}
This variational relation and the constraint (\ref{w0}) establish now the characterization of the range condition (\ref{res}) with the mixed variational problem (\ref{BB1}).  This characterization, along with (\ref{w2}) and (\ref{up56}), establishes Theorem \ref{well}$(ii)$.

By way of establishing wellposedness of (\ref{BB1}):  The bilinear forms $a_{\lambda }(\cdot ,\cdot )$ and $b(\cdot, \cdot)$ are readily
seen to be continuous.  In addition $a_\lambda$ is $H_{0}^{2}(\Omega )$-elliptic. The existence of a unique pair $[w_{1},\tilde{c%
}]\in H_{0}^{2}(\Omega )\times \mathbb{R}$ which solves (\ref{BB1}) will follow the Babu\v{s}ka-Brezzi Theorem if we establish the following
\textquotedblleft inf-sup\textquotedblright condition, for some positive constant $\beta $:
\begin{equation}
\sup_{\phi \in H_{0}^{2}(\Omega )}\frac{b(\phi ,r)}{\Vert \phi \Vert
_{H_{0}^{2}(\Omega )}}\geq \beta |r|,\hspace{0.1in}\forall r\in \mathbb{R}.
\label{infsup}
\end{equation}%
To this end, consider the function $\xi \in H^{4}(\Omega ) \cap H_{0}^{2}(\Omega )$ which solves
\begin{equation}
\Delta ^{2}\xi =1\text{ in }\Omega ;\,\,\left. \xi \right\vert _{\partial \Omega }=\left. \frac{\partial \xi }{\partial \nu }\right\vert
_{\partial \Omega }=0\text{ on }\partial \Omega . \label{xi}
\end{equation}%
By Green's Formula we then have
\begin{equation}
\int_{\Omega }\xi \cdot 1\,d\Omega =\int_{\Omega }\xi \Delta ^{2}\xi
\,d\Omega =\int_{\Omega }\Delta \xi \Delta \xi \,d\Omega .  \label{wbar}
\end{equation}%
Therewith, for given scalar $r\in \mathbb{R}$, let $\eta \equiv -\text{sgn}(r)\xi $%
. Then
\begin{align}
\sup_{\phi \in H_{0}^{2}(\Omega )}\frac{b(\phi ,r)}{\Vert \phi \Vert
_{H_{0}^{2}(\Omega )}}& \geq \frac{-r\int_{\Omega }\eta \,d\Omega }{\Vert
\eta \Vert _{H_{0}^{2}(\Omega )}} \label{infsupb}\\
& =\frac{|r|\int_{\Omega }\xi \,d\Omega }{\Vert \xi \Vert _{H_{0}^{2}(\Omega
)}} \notag \\
& =\Vert \xi \Vert _{H_{0}^{2}(\Omega )}|r|, \notag
\end{align}%
after using (\ref{wbar}). This gives (\ref{infsup}), with inf-sup constant $%
\beta =\Vert \xi \Vert _{H_{0}^{2}(\Omega )}$. The existence of a unique pair $[w_{1},\tilde{c}]$ which solves (\ref{BB1}) now follows from
Theorem \ref{thm:BB}.

Note in particular that
\begin{equation}
w_{1}\in \frac{L^{2}(\Omega )}{\mathbb{R}}\,,  \label{mean}
\end{equation}%
from the second equation of (\ref{BB1}). In turn, we recover $w_{2},u,$ and $%
p$ via
\begin{equation}
\begin{cases}
w_{2}\equiv \lambda w_{1}-w_{1}^{\ast }\in H_{0}^{2}(\Omega )\cap \frac{%
L^{2}(\Omega )}{\mathbb{R}}, \\
u\equiv \tilde{f}(\lambda w_{1}-w_{1}^{\ast })+\tilde{\mu}(u^{\ast })\in
\mathbf{H}^{2}(\mathcal{O})\cap \mathcal{H}_{\text{fluid}}, \\
p\equiv \tilde{\pi}(\lambda w_{1}-w_{1}^{\ast })+\tilde{q}(u^{\ast })+\tilde{%
c}\in H^{1}(\mathcal{O}).%
\end{cases}
\label{three}
\end{equation}

\indent From (\ref{sol1}) and (\ref{sol2}),
\begin{equation}
\text{the variables }u\text{ and }p\text{ solve the Stokes system (\ref{flu1}%
).}  \label{four}
\end{equation}%
Moreover, from (\ref{var1}) and (\ref{sol1}) we have
\begin{align*}
& \lambda ^{2}(P_{\rho }^{1/2}w_{1},P_{\rho }^{1/2}\phi )_{\Omega }+(\Delta
w_{1},\Delta \phi )_{\Omega }+\lambda (u,\tilde{f}(\phi ))_{\mathcal{O}%
}+(\nabla u,\nabla \tilde{f}(\phi ))_{\mathcal{O}}-(p,\mathrm{div}(\tilde{f}%
(\phi )))_{\mathcal{O}} \\
& =\hspace{0.1in}(u^{\ast },\tilde{f}(\phi ))_{\mathcal{O}}+\left( P_{\rho
}(\lambda w_{1}^{\ast }+w_{2}^{\ast }),\phi \right) _{\Omega }\hspace{0.2in}%
\forall \phi \in H_{0}^{2}(\Omega ).
\end{align*}%
An integration by parts and passing of the adjoint of $P_{\rho }^{1/2}$ yields,
\begin{align*}
& \lambda ^{2}(P_{\rho }w_{1},\phi )_{\Omega }+(\Delta w_{1},\Delta \phi
)_{\Omega }+(\lambda u-\Delta u+\nabla p,\tilde{f}(\phi ))_{\mathcal{O}%
}+\left\langle \frac{\partial u}{\partial \nu }-p\nu ,\tilde{f}(\phi
)\right\rangle _{\partial \mathcal{O}} \\
& =\hspace{0.1in}(u^{\ast },\tilde{f}(\phi ))_{\mathcal{O}}+\left( P_{\rho
}(\lambda w_{1}^{\ast }+w_{2}^{\ast }),\phi \right) _{\Omega }\hspace{0.2in}%
\forall \phi \in H_{0}^{2}(\Omega ).
\end{align*}%
As variables $u$ and $p$ solve the Stokes system (\ref{flu1}), we thus
attain the relation
\begin{equation*}
\lambda ^{2}(P_{\rho }w_{1},\phi )_{\Omega }+(\Delta w_{1},\Delta \phi
)_{\Omega }-(\left. p\right\vert _{\Omega },\phi )_{\Omega }=\left( P_{\rho
}(\lambda w_{1}^{\ast }+w_{2}^{\ast }),\phi \right) _{\Omega },\hspace{0.2in}%
\forall \phi \in H_{0}^{2}(\Omega ).
\end{equation*}%
(We have also implicitly used (\ref{sol1}) and the remark after (\ref{wk1}%
).) In particular, this holds true for $\phi \in \mathcal{D}(\Omega )$. Thus
we have the distributional relation
\begin{equation*}
(\lambda ^{2}P_{\rho }w_{1}+\Delta ^{2}w_{1}-\left. p\right\vert _{\Omega
}-P_{\rho }[\lambda w_{1}^{\ast }+w_{2}^{\ast }],\phi )_{\Omega }=0\hspace{%
0.2in}\forall \phi \in \mathcal{D}(\Omega ),
\end{equation*}%
and so we infer that
\begin{equation}
\begin{array}{c}
w_{1}\text{ satisfies (\ref{str1}).}%
\end{array}
\label{str2}
\end{equation}%
Subsequently, we infer by elliptic theory that, as required by the
definition of the fluid-structure operator $\mathcal{A}_{\rho }:D(\mathcal{A}%
_{\rho })\subset \mathbf{H}_{\rho }\rightarrow \mathbf{H}_{\rho }$,%
\begin{equation}
w_{1}\in \mathcal{S}_{\rho },  \label{e_reg}
\end{equation}%
where the (displacement) space is as given in (\ref{S}). Finally,
because $u$ and $u^{\ast }\in \mathcal{H}_{\mathrm{fluid}}$, we have \textit{%
a fortiori} from (\ref{flu1}),%
\begin{equation}
\Delta p=0\text{ in }\mathcal{O}\,\text{ and }\frac{\partial p}{\partial \nu
}\Big|_{S}=\Delta u\cdot \nu \big|_{S}\text{ on }S.  \label{13}
\end{equation}%
Moreover, from (\ref{three}) and (\ref{str1}), $\lambda w_{2}+P_{\rho
}^{-1}\Delta ^{2}w_{1}-P_{\rho }^{-1}p\big|_{\Omega }=w_{2}^{\ast }$ which
implies that in $\Omega $, (since $[w_{1}^{\ast },w_{2}^{\ast },u^{\ast
}]\in \mathbf{H}_{\rho }$)
\begin{align}
P_{\rho }^{-1}\Delta ^{2}w_{1}& =P_{\rho }^{-1}p\big|_{\Omega }-\lambda
w_{2}+w_{2}^{\ast }  \notag \\
& =P_{\rho }^{-1}p\big|_{\Omega }-\lambda u\cdot \nu \big|_{\Omega }+u^{\ast
}\cdot \nu \big|_{\Omega }  \notag \\
& =P_{\rho }^{-1}p\big|_{\Omega }-\Delta u\cdot \nu \big|_{\Omega }+\nabla
p\cdot \nu \big|_{\Omega }  \label{14}
\end{align}%
where in the last equality, (\ref{flu1}) was again invoked. Thus, from (\ref%
{13}) and (\ref{14}), we have that the pressure variable $p$ we have obtained by Theorem \ref{thm:BB} solves
\begin{equation*}
\begin{cases}
\Delta p=0 & \text{in }\mathcal{O}, \\
\frac{\partial p}{\partial \nu }+P_{\rho }^{-1}p=P_{\rho }^{-1}\Delta
^{2}w_{1}+\Delta u^{3}\big|_{\Omega } & \text{in }\Omega , \\
\frac{\partial p}{\partial \nu }=\Delta u\cdot \nu \big|_{S} & \text{on }S.%
\end{cases}%
\end{equation*}%
As such,%
\begin{equation}
p=G_{\rho ,1}(w_{1})+G_{\rho ,2}(u),  \label{press}
\end{equation}%
where $G_{\rho ,i}$ are as given in (\ref{G1}) and (\ref{G2}). (Note that we
are implicitly using the critical regularity
\begin{equation}
\Delta u\cdot \nu \big|_{\partial \mathcal{O}}\in H^{-\frac{1}{2}}(\partial \mathcal{O%
}),  \label{trace}
\end{equation}%
from (\ref{CritReg}).)

\indent From (\ref{mean}), (\ref{three}), (\ref{four}), (\ref{str2}), (\ref{e_reg}),
and (\ref{press}), we have that the constructed variables $%
[w_{1},w_{2},u]$ belong to $D(\mathcal{A})$ and solve the
resolvent equation (\ref{res}). This concludes the proof of Theorem \ref%
{well}. \ \ \ $\square $

\section{A Numerical Analysis of the Fluid-Structure Dynamics}

The objective of this section is to demonstrate how the maximality argument which was given in Section \ref{Maximality} can be utilized to
approximate solutions to the fluid-structure interactive PDE under present consideration.  In particular, the numerical method outlined here
solves the static problem resulting from the resolvent equations (\ref{w2}) - (\ref{end}), but this approach in principle can be modified to
solve the time dependent problem in the same way as in \cite{dvorak} (see p. 276). We will outline here a certain numerical implementation of the
finite element method (FEM) and provide convergence results for the approximation with respect to ``mesh parameter" $h$.  Finally, we will
provide an explicit model problem as a numerical example.

\subsection{Finite Element Formulation}

In what follows, the three dimensional body $\mathcal{O}$ will be taken to be a polyhedron.  Given a positive (and small) parameter $h$ of
discretization, we let $\{e_\ell\}_{\ell = 1}^{N_h}$ be an FEM ``triangulation" of $\mathcal{O}$, where each element $e_\ell$
is a tetrahedron (and so, among other properties, $\displaystyle \bigcup_{\ell=1}^{N_h} e_\ell = \mathcal{O}$, see \cite{Axelsson} and Figure
\ref{GMSH} below).

\begin{enumerate}
\item[(A)]
Relative to the ``triangulation" of $\mathcal{O}$, $V_h$ will denote the classic $\mathbf{H}^1$-conforming FEM finite dimensional subspace such
that
\begin{equation}
[\mathbb{P}_2]^3 \subset V_h \subset \mathbf{H}^1_0(\mathcal{O}), \,\,\, V_h \not\subset \mathbf{H}^2(\mathcal{O}); \,\,\, V_h \subset
[\mathcal{C}(\bar{\Omega})]^3, \,\,\, V_h \not\subset [\mathcal{C}^1(\bar{\Omega})]^3. \label{Vh}
\end{equation}
(See \cite{Axelsson}.)  Subsequently to handle the inhomogeneity we specify the set
\begin{equation}
\tilde{V}_h = \Big\{\mu_h + \gamma^+_0(\xi) \in \mathbf{H}^1(\mathcal{O}): \mu_h \in V_h \text{ and } \gamma_0^+(\xi) \Big|_{\partial \mathcal{O}} = \begin{cases} \vec{0}, & \text{ on } S \\ [0,0,\xi], & \text{ for } \xi \in H^2(\Omega),\end{cases} \Big\}.
\end{equation}
\item[(B)]
In addition $\Pi_h$ will denote the $L^2$-FEM finite dimensional subspace for the pressure variable defined by
\begin{equation}
\Pi_h = \big\{ q_h \in \frac{L^2(\mathcal{O})}{\mathbb{R}} \cap \mathcal{C}(\bar{\mathcal{O}}): \,\, \forall \ell = 1,. . . , N_h; \,\, q_h
\big|_{e_\ell} \in \mathbb{P}^1 \big\} \label{Pih}
\end{equation}
(see \cite{Ciarlet} and \cite{Axelsson}). Moreover, we let $\{\tilde{e}_\ell\}_{\ell =1}^{\tilde{N}_h}$ be a FEM triangulation of the two
dimensional polygonal region $\Omega$, where each element $\tilde{e}_\ell$ is a triangle.
\item[(C)]
Similarly, $X_h$ will denote a conforming FEM subspace such that
\begin{equation}
\mathbb{P}_3 \subset X_h \subset H^2_0(\Omega), \,\,\, X_h \not\subset H^3(\Omega); \,\,\, X_h \subset \mathcal{C}^1(\bar{\Omega}), \,\,\, X_h
\not\subset \mathcal{C}^2(\bar{\Omega}) \label{Xh}
\end{equation}
(see e.g. \cite{Ciarlet} and \cite{Solin} for details of the explicit construction of these piecewise polynomials.  As such, the basis functions
which generate $X_h$ are ``conforming", relative to fourth-order boundary value problems.)
\end{enumerate}

For the spaces $V_h$, $\Pi_h$ and $X_h$ described above we will have need of the following discrete estimates relative to mesh parameter $h$:
\begin{enumerate}
\item[(A$'$)] In regard to the $\mathbf{H}^1(\mathcal{O})$-conforming FEM space $V_h$ in (\ref{Vh}) we have the following estimate:
For $\mu \in \mathbf{H}^2(\mathcal{O}) \cap \mathbf{H}^1_0(\mathcal{O}),$
\begin{equation}
\min_{\mu_h \in V_h} \|\mu - \mu_h\|_{\mathbf{H}^1_0(\mathcal{O})} \leq Ch |\mu|_{2,\mathcal{O}}. \label{muEstimate}
\end{equation}
(See Theorem 5.6, p. 224, of \cite{Axelsson}.)
\item[(B$'$)] Similarly, in regard to the finite dimensional space $\Pi_h$, we have the discrete estimate:
For $q \in H^1(\mathcal{O}) \cap \frac{L^2(\mathcal{O})}{\mathbb{R}},$
\begin{equation}
\min_{q_h \in \Pi_h} \|q - q_h\|_{L^2(\mathcal{O})} \leq Ch \|q\|_{H^1(\mathcal{O})}. \label{qEstimate}
\end{equation}
(See e.g., Corollary 1.128, p. 70,  of \cite{Ern}.)
\item[(C$'$)] Finally, with regard to the FEM space $X_h$ in (\ref{Xh}), we have the following discrete estimates:
\begin{enumerate}
\item[(i)] For $\psi \in H^4(\Omega) \cap H^2_0(\Omega)$,
\begin{equation}
\min_{\psi_h \in X_h} \|\psi - \psi_h\|_{H^2_0(\Omega)} \leq C h^2 |\psi|_{4, \Omega}, \label{XhEstimate1}
\end{equation}
\item[(ii)] For $\psi \in H^3(\Omega) \cap H^2_0(\Omega)$,
\begin{equation}
\min_{\psi_h \in X_h} \|\psi - \psi_h\|_{H^2_0(\Omega)} \leq C h |\psi|_{3, \Omega}. \label{XhEstimate2}
\end{equation}
(See estimate (5.82), p. 225, of \cite{Axelsson}.)
\end{enumerate}
\end{enumerate}

The goal here is to find a finite dimensional approximation $[w_{1h}, w_{2h}, u_h] \in X_h \times X_h \times \tilde{V}_h$ to the solution $[w_1,
w_2, u] \in D(\mathcal{A}_\rho)$ of (\ref{res}), as well as an approximation $p_h$ of the associated fluid pressure $p$.  We shall see that
these particular FEM subspaces are chosen with a view of satisfying the (discrete) Babu\v{s}ka-Brezzi condition relative to a mixed variational
formulation, a formulation which is wholly analogous to that in (\ref{BB1}) for the static fluid-structure PDE system (\ref{w2})-(\ref{end}). We
further note that, by way of satisfying said inf-sup condition, it is indispensable that the structural component space $X_h$ be
$H^2$-conforming (see (\ref{Xh})). In addition, this mixed variational formulation for the coupled problem (\ref{w2})-(\ref{end}), like the
mixed method for uncoupled Stokes or Navier-Stokes flow, allows for the implementation of approximating fluid basis functions (in $V_h$) which
are \underline{not} divergence free (see \cite{Fortin}).

In line with the maximality argument of Section \ref{Maximality}, the initial task in the present finite dimensional setting is to numerically
resolve the structural solution component of the PDE system (\ref{w2})-(\ref{end}).  Namely, with reference to the bilinear and linear
functionals $a_\lambda(\cdot, \cdot), b(\cdot, \cdot)$ and $\mathbb{F}(\cdot)$ of (\ref{abF}), the present discrete problem is to find
$[w_{1h},\tilde{c}_h] \in X_h \times \mathbb{R}$ which solve:
\begin{equation}
\begin{cases}
a_\lambda(w_{1h}, \psi_h) + b(\psi_h, \tilde{c}_h) = \mathbb{F}(\psi_h), & \forall \psi_h \in X_h,\\
b(w_{1h}, r) = 0, & \forall r \in \mathbb{R}.
\end{cases}\label{BB1h}
\end{equation}

Assuming this variational problem can be solved uniquely, $w_{2h}$ is immediately resolved via the relation
\begin{equation}
w_{2h} = \lambda w_{1h} - w_1^* \label{w2h}
\end{equation}
(cf. (\ref{w2})).  Subsequently we can recover fluid and pressure approximations $u_h$ and $p_h$ from the discrete solution pair $[w_{1h},
\tilde{c}_h] \in X_h \times \mathbb{R}$ of (\ref{BB1h}).  Indeed, to this end we will invoke the classic mixed variational formulation for
Stokes flow, so as to approximate the fluid maps $[\tilde{f}(\cdot), \tilde{\pi}(\cdot)]$ and $[\tilde{\mu}(\cdot), \tilde{q}(\cdot)]$ of
(\ref{sol1}) and (\ref{sol2}), respectively.  (See \cite{Fortin}.) Let bilinear forms $\tilde{\mathbf{a}}_\lambda (\cdot, \cdot):
\mathbf{H}^1_0(\mathcal{O}) \times \mathbf{H}^1_0(\mathcal{O}) \to \mathbb{R}$ and $\tilde{\mathbf{b}}(\cdot, \cdot) :
\mathbf{H}^1_0(\mathcal{O}) \times \mathbf{L}^2(\mathcal{O})/\mathbb{R} \to \mathbb{R}$ be defined respectively as follows:
\begin{alignat}{3}
\tilde{\mathbf{a}}_\lambda(\mu, \varphi) &= \lambda (\mu, \varphi)_\mathcal{O} + (\nabla \mu, \nabla \varphi)_\mathcal{O}, &\quad& &&\forall
\mu, \varphi \in \mathbf{H}^1_0(\mathcal{O}); \label{atilde}\\
\tilde{\mathbf{b}}(\mu, q) &= -(\text{div}(\mu), q)_\mathcal{O}, &\quad& &&\forall \mu \in \mathbf{H}^1_0(\mathcal{O}), q \in
\frac{L^2(\mathcal{O})}{\mathbb{R}}. \label{btilde}
\end{alignat}

Moreover, we define the standard Sobolev trace map $\gamma_0: \mathbf{H}^k(\mathcal{O}) \to \mathbf{H}^{k-1/2}(\partial \mathcal{O})$, for
$k=1,2,3, \hdots$. That is for $f \in [\mathcal{C}^\infty(\bar{\mathcal{O}})]^3$,
\begin{equation*}
\gamma_0(f) = f\big|_{\partial \mathcal{O}}.
\end{equation*}
Since $\gamma_0(\cdot)$ is continuous and surjective, then for any $\phi \in H^{k-1/2}(\Omega)$ we have the existence and uniqueness of an element
in $\mathbf{H}^k(\mathcal{O})$, denoted here as $\gamma_0^+(\phi)$, which satisfies
\begin{equation}
\gamma_0 \gamma_0^+(\phi) = \begin{cases} \vec{0} & \text{ on } S,\\ [0, 0, \phi] & \text{ on } \Omega.
\end{cases} \label{gamma+}
\end{equation}

Therewith, the classic mixed FEM for (\ref{sol1}) is given as follows: With subspaces $V_h$ and $\Pi_h$ as given in (\ref{Vh}) and (\ref{Pih})
respectively, and given $\phi \in H^{1/2}(\Omega)$, find the unique pair $[\tilde{f}_{0h}(\phi), \tilde{\pi}_h(\phi)] \in V_h \times \Pi_h$ such
that
\begin{alignat}{3}
\tilde{\mathbf{a}}_\lambda (\tilde{f}_{0h},\varphi_h) + \tilde{\mathbf{b}}(\varphi_h, \tilde{\pi}_{h}) &= - \tilde{\mathbf{a}}_\lambda(
\gamma_0^+(\phi), \varphi_h) &\quad & && \forall \varphi_h \in V_h, \label{ah1}\\
\tilde{\mathbf{b}}(\tilde{f}_{0h}, \varrho_h) &=-\left[\frac{\int_\Omega \phi \, d\Omega}{\text{meas}(\mathcal{O})}\right] \int_\mathcal{O}
\varrho_h \, d\mathcal{O} - \tilde{\mathbf{b}}(\gamma_0^+(\phi),\varrho_h) &\quad & && \forall \varrho_h \in \Pi_h. \label{bh1}
\end{alignat}

Likewise, the classic mixed FEM for (\ref{sol2}) is given as follows: For given $u^* \in H^{-1}(\mathcal{O})$, find the unique pair
$[\tilde{\mu}_h(u^*), \tilde{q}_h(u^*)] \in V_h \times \Pi_h$ such that
\begin{alignat}{3}
\tilde{\mathbf{a}}_\lambda (\tilde{\mu}_{h},\varphi_h) + \tilde{\mathbf{b}}(\varphi_h, \tilde{q}_{h}) &= (u^*, \varphi_h)_\mathcal{O} &\quad &
&&
\forall \varphi_h \in V_h; \label{ah2}\\
\tilde{\mathbf{b}}(\tilde{\mu}_{h}, \varrho_h) &= 0 &\quad & && \forall \varrho_h \in \Pi_h. \label{bh2}
\end{alignat}

By the Babu\v{s}ka-Brezzi Theorem, the two discrete variational formulations (\ref{ah1})-(\ref{bh1}) and (\ref{ah2})-(\ref{bh2}) are well-posed;
see \cite{Fortin}. (In particular, with the so-called Taylor-Hood formulation in place - i.e, fluid approximation space $V_h$ consists of
piecewise quadratic functions, and pressure approximation space $\Pi_h$ consists of piecewise linear functions - then the aforesaid inf-sup
condition is satisfied \emph{uniformly} in parameter $h$.)

With the approximating solution maps (\ref{ah1})-(\ref{bh2}) in place and assuming the structural component approximation $[w_{1h},
\tilde{c}_h]$ is known, we then set
\begin{align}
u_h &= \tilde{f}_{0h}(\lambda w_{1h} - w_1^*)+ \gamma_0^+(\lambda w_{1h} - w_1^*) + \tilde{\mu}_{h}(u^*) ; \label{uh}\\
p_h &= \tilde{\pi}_{h}(\lambda w_{1h} - w_1^*) + \tilde{q}_{h}(u^*) + \tilde{c}_h,\label{ph}
\end{align}
(c.f. (\ref{up}).)

Now in regard to the variational problem in (\ref{BB1h}), one will in fact have unique solvability of this discrete problem, via the
Babu\v{s}ka-Brezzi Theorem, provided that the following inf-sup condition is satisfied:
\begin{equation}
\sup_{\phi_h \in X_h} \frac{b(\phi_h, r)}{\|\phi_h\|_{H^2_0(\Omega)}} \geq \beta_h |r|, \, \forall r \in \mathbb{R}. \label{infsup1}
\end{equation}
But what is more, in order to ensure stability and ultimately convergence of the numerical solutions obtained by our particular FEM, it is
indispensable that the ``discrete" inf-sup condition (\ref{infsup1}) be \emph{uniform} of parameter $h>0$ (at least for $h$ small enough).

In fact we have the following result:
\begin{lemma}
\label{lemma:infsup} Let the bilinear form $b(\cdot, \cdot): H_0^2 (\Omega) \times \mathbb{R} \to \mathbb{R}$ be as defined in (\ref{abF}). Then
for parameter $h >0$ small enough one has the ``inf-sup" estimate
\begin{equation}
\sup_{\phi_h \in X_h} \frac{b(\phi_h, r)}{\|\phi_h\|_{H^2_0(\Omega)}} \geq C |r|, \, \forall r \in \mathbb{R}. \label{infsup2}
\end{equation}
where $C = \|\xi\|_{H^2_0(\Omega)} - \epsilon$, and $\xi$ is the solution of the boundary value problem (\ref{xi}).  Here, $\epsilon>0$ can be
taken arbitrarily small.
\end{lemma}
\textbf{Proof of Lemma \ref{lemma:infsup}.}  We resurrect the elliptic variable $\xi \in  H^4(\Omega) \cap H^2_0(\Omega)$ from the earlier
maximality argument. Namely, $\xi$ solves
\begin{equation*}
\Delta^2 \xi = 1 \text{ in } \Omega; \,\,\,\,\, \xi\big|_{\partial \Omega} = \frac{\partial \xi}{\partial \nu} \Big|_{\partial \Omega} = 0
\text{ on }
\partial \Omega.
\end{equation*}
Then by Green's First Identity we have that $\xi$ solves the following variational problem for all $\psi \in H^2_0(\Omega)$:
\begin{align}
(1, \psi)_\Omega &= (\Delta^2 \xi, \psi)_\Omega \notag\\
&= (\Delta \xi, \Delta \psi)_\Omega, \text{ for all } \psi \in H^2_0(\Omega). \label{xiIdentity}
\end{align}
Let now $\xi_h \in X_h$ denote the ``energy projection" of $\xi$ on $X_h$.  That is, $\xi_h$ satisfies the following discrete variational
problem:
\begin{equation}
(\Delta \xi_h, \Delta \psi_h)_\Omega = (1, \psi_h)_\Omega \,\,\, \forall \, \psi_h \in X_h. \label{xihIdentity}
\end{equation}
The existence and uniqueness of the discrete solution $\xi_h \in X_h$ follows from the Lax-Milgram Theorem, see e.g., \cite{Axelsson},
\cite{Ciarlet}.  Applying the discrete estimate (\ref{XhEstimate1}) to the respective variational problems (\ref{xiIdentity}) and
(\ref{xihIdentity}), we then have
\begin{equation}
|\Delta(\xi - \xi_h)|_\Omega \leq Ch^2. \label{xiEstimate}
\end{equation}
With these ingredients, $\xi$ and $\xi_h$, we then have
\begin{align}
\sup_{\phi_h \in X_h} \frac{b(\phi_h, r)}{\|\phi_h\|_{H^2_0(\Omega)}} &\geq \frac{-r \int_\Omega [-\text{sgn}(r)]\xi_h
\, d\Omega}{\|\xi_h\|_{H^2_0(\Omega)}} \notag\\
&= \frac{|r| \int_\Omega \xi_h\cdot 1 \, d\Omega}{\|\xi_h\|_{H^2_0(\Omega)}} \notag\\
&= |r| \big\|\xi_h\big\|_{H^2_0(\Omega)},
\end{align}
after using (\ref{xihIdentity}) above.  Continuing, we then have
\begin{align}
\sup_{\phi_h \in X_h} \frac{b(\phi_h, r)}{\|\phi_h\|_{H^2_0(\Omega)}} &\geq |r| \big\| \xi - (\xi - \xi_h)\big\|_{H^2_0(\Omega)} \notag\\
&\geq |r| \big[\|\xi\|_{H^2_0(\Omega)} - \|\xi - \xi_h\|_{H^2_0(\Omega)}\big] \notag \\
&\geq |r| \big[\|\xi\|_{H^2_0(\Omega)} - Ch^2\big],
\end{align}
after using the estimate (\ref{xiEstimate}).  Taking step size parameter
\begin{equation}
h < \sqrt{\frac{\epsilon}{C}}
\end{equation}
now completes the proof. \ \ \ $\square $

\subsection{Error Estimates for the Finite Element Problem}
In what follows we will have need of the following result in \cite{Ern}, which will not be stated here in its full generality (see \cite{Ern},
Lemma 2.44, p. 104).
\begin{lemma}
\label{lemma:2.44} With reference to the quantities in Theorem \ref{thm:BB} above, let $\Sigma_h$ be a subspace of $\Sigma$, and let $M_h$ be a
subspace of $M$.  Suppose further that bilinear form $a: \Sigma \times \Sigma \to \mathbb{R}$ is $\Sigma$-elliptic; that is, $\exists \, \alpha
>0$ such that
\begin{equation}
a(\sigma, \sigma) \geq \alpha \|\sigma \|_\Sigma^2 \,\,\,\, \forall \sigma \in \Sigma . \label{SigmaElliptic}
\end{equation}
Also, assume that the following ``discrete inf-sup" condition is satisfied: $\exists \, \beta_h >0$ such that
\begin{equation}
\inf_{q_h \in M_h} \sup_{\tau_h \in \Sigma_h} \frac{b(\tau_h, q_h)}{\|\tau_h\|_\Sigma \|q_h\|_{M}} \geq \beta_h, \label{ErnInfSup}
\end{equation}
where $\beta_h>0$ may depend upon subspaces $\Sigma_h$ and $M_h$.  Let moreover $(\sigma_h, p_h) \in \Sigma_h \times M_h$ solve the following
(approximating) variational problem:
\begin{equation}
\begin{cases} a(\sigma_h, \tau_h) + b(\tau_h, p_h) = (\kappa, \tau_h) &  \forall \tau_h \in \Sigma_h, \\
b(\sigma_h, q_h) = (\ell, q_h) & \forall q_h \in M_h.
\end{cases}\label{ErnVariational1}
\end{equation}
(Note that the existence and uniqueness of the solution pair $(\sigma_h, p_h)$ follows from Theorem \ref{thm:BB}, in view of (\ref{SigmaElliptic}) and (\ref{ErnInfSup}).) Then one has the following error estimates:
\begin{align}
\|\sigma - \sigma_h\|_\Sigma &\leq c_{1h} \inf_{\varsigma_h \in \Sigma_h} \|\sigma - \varsigma_h\|_\Sigma, + c_{2h} \inf_{q_h \in M_h} \|p -
q_h\|_M \label{ErnVariational2} \\
\|p - p_h\|_M &\leq c_{3h} \inf_{\varsigma_h \in \Sigma_h} \|\sigma - \varsigma_h\|_\Sigma + c_{4h} \inf_{q_h \in M_h} \|p - q_h\|_M,
\label{ErnVariational3}
\end{align}
with $c_{1h} = (1+\frac{\|a\|}{\alpha_h})(1+\frac{\|b\|}{\beta_h})$, $c_{2h} = \frac{\|b\|}{\alpha_h}$; moreover if $M = M_h$ one can take $c_{2h} = 0$, $c_{3h} = c_{1h}\frac{\|a\|}{\beta_h}$,
and $c_{4h} = 1 + \frac{\|b\|}{\beta_h} + c_{2h}\frac{\|a\|}{\beta_h}$.
\end{lemma}

Concerning the efficacy of our FEM for numerically approximating the fluid-structure system (\ref{w2})-(\ref{end}), we have the following:

\begin{theorem}\label{thm:estimates}
Let $h>0$ be the parameter of discretization which gives rise to the FEM subspaces $V_h$, $\Pi_h$, and $X_h$ of (\ref{Vh}), (\ref{Pih}), and
(\ref{Xh}), respectively. With respect to the solution variables $[w_1, w_2, u, p] \in D(\mathcal{A}_\rho)\times H^1(\mathcal{O})$ of
(\ref{w2})-(\ref{end}) and their FEM approximations $[w_{1h}, w_{2h}, u_h, p_h]$, as given by (\ref{BB1h}) - (\ref{w2h}) and
(\ref{uh})-(\ref{ph}), we have the following rates of convergence:
\begin{enumerate}
\item[(i)]  (a) If $\rho = 0$,
\begin{equation} \|w_1 - w_{1h}\|_{H^2_0(\Omega)} \leq C_\lambda h^2 \big\|[w_1^*, w_2^*, u^*]\big\|_{\mathbf{H}_0}. \end{equation}
(b) If $\rho > 0$,
\begin{equation} \|w_1 - w_{1h}\|_{H^2_0(\Omega)} \leq C_\lambda h \big\|[w_1^*, w_2^*, u^*]\big\|_{\mathbf{H}_\rho}. \end{equation}
\item[(ii)](a$'$) If $\rho = 0$,
\begin{equation} \|w_2 - w_{2h}\|_{H^2_0(\Omega)} \leq C_\lambda h^2 \big\|[w_1^*, w_2^*, u^*]\big\|_{\mathbf{H}_0}. \end{equation}
(b$'$) If $\rho > 0$,
\begin{equation} \|w_2 - w_{2h}\|_{H^2_0(\Omega)} \leq C_\lambda h \big\|[w_1^*, w_2^*, u^*]\big\|_{\mathbf{H}_\rho}. \end{equation}
\item[(iii)] For $\rho \geq 0$ \begin{equation}\|u-u_h\|_{\mathbf{H}^1(\mathcal{O})} \leq C_\lambda h \big\|[w_1^*, w_2^*,
u^*]\big\|_{\mathbf{H}_\rho}.\end{equation}
\item[(iv)] For $\rho \geq 0$ \begin{equation}\|p-p_h\|_{L^2(\mathcal{O})} \leq C_\lambda h \big\|[w_1^*, w_2^*, u^*]\big\|_{\mathbf{H}_\rho} \label{ivpressure} \end{equation}
\end{enumerate}
\end{theorem}

\textbf{Proof of Theorem \ref{thm:estimates}.}  We first establish parts $(i)$ and $(ii)$ together.  Here we will combine Lemma \ref{lemma:2.44}
with the bilinear forms in (\ref{abF}).  We take in Lemma \ref{lemma:2.44}
\begin{equation*}
a(\cdot, \cdot) \equiv a_\lambda(\cdot, \cdot): H^2_0(\Omega) \times H^2_0(\Omega) \to \mathbb{R},
\end{equation*}
so $\alpha = 1$ and $\|a\|= \|a_\lambda\|_{\mathcal{L}\left([H^2_0(\Omega)]^2, \mathbb{R}\right)} \leq C_\lambda$.  Moreover, we as before set
\begin{equation*}
b(\cdot,\cdot): H^2_0(\Omega) \times \mathbb{R} \to \mathbb{R} \text{ as in (\ref{abF})}
\end{equation*}
so $\|b\|_{\mathcal{L}\left(H^2_0(\Omega)\times \mathbb{R}, \mathbb{R}\right)} \leq C_\lambda$.

In addition, by Lemma \ref{lemma:infsup}, we can take
\begin{equation*}
\beta_h \equiv C, \,\,\, \forall h>0.
\end{equation*}
Subsequently, with reference to the variational system (\ref{BB1}) and the approximating system (\ref{BB1h}) we have from
(\ref{ErnVariational2})
\begin{equation}
\|w_1-w_{1h}\|_{H^2_0(\Omega)} \leq C_\lambda \inf_{\psi_h \in X_h} \|w_1 - \psi_h\|_{H^2_0(\Omega)};
\end{equation} (note that that second term from the right hand side of (\ref{ErnVariational2}) is zero in this case because $M_h = M
= \mathbb{R}$ in this case).  Appealing now to estimates (\ref{XhEstimate1}), (\ref{XhEstimate2}), and Theorem \ref{well}$(i)$ we have the following
error estimates:
\begin{enumerate}
\item[$(a)$] If $\rho =0$,
\begin{align}
\|w_1 - w_{1h}\|_{H^2_0(\Omega)} &\leq C_\lambda h^2 |w_1|_{4,\Omega}, \notag \\
&\leq C_\lambda h^2 \big\| [w_1^*, w_2^*, u^*]\big\|_{\mathbf{H}_0}. \label{w1Est1}
\end{align}
\item[$(b)$] If $\rho>0$, then
\begin{align}
\|w_1 - w_{1h}\|_{H^2_0(\Omega)} &\leq C_\lambda h |w_1|_{3,\Omega}, \notag \\
&\leq C_\lambda h \big\| [w_1^*, w_2^*, u^*]\big\|_{\mathbf{H}_\rho}. \label{w1Est2}
\end{align}
\end{enumerate}
This establishes Theorem \ref{thm:estimates}$(i)$.  In view of (\ref{w2h}), Theorem \ref{thm:estimates}$(ii)$ follows directly.

To achieve the estimates for the fluid variables, we first note that the error in the constant component of the pressure, given by the
variational system (\ref{BB1}), satisfies:
\begin{enumerate}
\item[$(a')$] If $\rho = 0$, we have upon combining (\ref{ErnVariational3}) and (\ref{w1Est1})
\begin{equation} |c - \tilde{c}_h| \leq C_\lambda h^2 \big\|[w_1^*, w_2^*, u^*]\big\|_{\mathbf{H}_0} .\label{cEst3}
\end{equation}
(Implicity we are taking in (\ref{ErnVariational3}) [finite dimensional] $\mathbb{R} = M_h$.)
\item[$(b')$] If $\rho >0$, we have upon combining (\ref{ErnVariational3}) and (\ref{w1Est2})
\begin{equation} |c - \tilde{c}_h| \leq C_\lambda h \big\|[w_1^*, w_2^*, u^*]\big\|_{\mathbf{H}_\rho} .\label{cEst4}
\end{equation}
\end{enumerate}
(After again noting as above that in this case $M_h = M = \mathbb{R}$.)

Now we will again invoke Lemma \ref{lemma:2.44}, with respect to the Stokes component (\ref{ustar})-(\ref{end}), with therein
\begin{equation}
a(\cdot, \cdot) \equiv \tilde{\mathbf{a}}_\lambda(\cdot, \cdot) : \mathbf{H}^1_0(\mathcal{O}) \times \mathbf{H}^1_0(\mathcal{O}) \to \mathbb{R}
\text{ as given in (\ref{atilde})}.
\end{equation}
Consequently we can take $\alpha \equiv 1$ and $\|a\| = \|\tilde{\mathbf{a}}_\lambda \|_{\mathcal{L}\left([\mathbf{H}_0^1(\mathcal{O})]^2,
\mathbb{R}\right)} \leq C_\lambda.$  Moreover we take
\begin{equation}
b(\cdot, \cdot) \equiv \tilde{\mathbf{b}}(\cdot, \cdot) : \mathbf{H}_0^1(\mathcal{O}) \times \frac{L^2(\mathcal{O})}{\mathbb{R}} \to \mathbb{R},
\text{ as given in (\ref{btilde}).}
\end{equation}
Then $\|b\| = \|\tilde{\mathbf{b}}\|_{\mathcal{L}\left(\mathbf{H}_0^1(\mathcal{O}) \times \frac{L^2(\mathcal{O})}{\mathbb{R}},
\mathbb{R}\right)} \leq C.$

In addition, since the respective fluid and pressure spaces $V_h$ and $\Pi_h$ are piecewise quadratic and piecewise linear - i.e. the so-called
Taylor-Hood formulation - then for $h>0$ small enough we can take
\begin{equation*}
\beta_h = \beta^*,
\end{equation*}
independent of small $h>0$. (See e.g., Lemma 4.23, p. 193 of \cite{Ern}.)

Thus with $u$ and $u_h$ as given in (\ref{up}) and (\ref{uh}) respectively, we then have
\begin{align}
\|u- u_h\|_{\mathbf{H}^1(\mathcal{O})} & \leq \|\tilde{f}(\lambda w_1 - w_1^*) - \tilde{f}_{0h}(\lambda w_{1h} - w_1^*) - \gamma_0^+(\lambda
w_{1h} - w_1^*)\|_{\mathbf{H}^1(\mathcal{O})} + \|\tilde{\mu}(u^*) - \tilde{\mu}_h(u^*)\|_{\mathbf{H}^1(\mathcal{O})} \notag \\
&\leq \|\tilde{f}(\lambda w_1 - w_1^*) - \tilde{f}_{0h}(\lambda w_{1h} - w_1^*) - \gamma_0^+(\lambda w_{1h} -
w_1^*)\|_{\mathbf{H}^1(\mathcal{O})} + Ch\big(|u|_{2,\mathcal{O}} + \|p\|_{H^1(\mathcal{O})}\big)\label{uhEst1}
\end{align}
after using estimate (\ref{ErnVariational2}) followed by (\ref{muEstimate}) and (\ref{qEstimate}).

Concerning the first term on the right hand side of (\ref{uhEst1}) we have further
\begin{align}
\|\tilde{f}&(\lambda w_1 - w_1^*) - \tilde{f}_{0h}(\lambda w_{1h} - w_1^*) - \gamma_0^+(\lambda w_{1h} - w_1^*)\|_{\mathbf{H}^1(\mathcal{O})}
\notag \\
& \leq \|\tilde{f}(\lambda w_1 - w_1^*) - \tilde{f}_{0h}(\lambda w_{1} - w_1^*) - \gamma_0^+(\lambda w_{1} -
w_1^*)\|_{\mathbf{H}^1(\mathcal{O})}
\notag \\
&\quad + \|\gamma_0^+(\lambda w_{1} - w_1^*) - \gamma_0^+(\lambda w_{1h} - w_1^*)\|_{\mathbf{H}^1(\mathcal{O})} + \|\tilde{f}_{0h}(\lambda w_{1}
- w_1^*) - \tilde{f}_{0h}(\lambda w_{1h} - w_1^*)\|_{\mathbf{H}^1(\mathcal{O})} \label{uhEst2}
\end{align}
We now estimate these terms one at at time.

Appealing again to the estimates in (\ref{muEstimate}) and (\ref{qEstimate}) via (\ref{ErnVariational2}) (as well as to the regularity given in
(\ref{map1})) we have
\begin{align}
\big\|[\tilde{f}(\lambda w_1 - w_1^*) - \gamma_0^+(\lambda w_{1} - w_1^*)] - \tilde{f}_{0h}(\lambda w_{1} - w_1^*)\big\|_{\mathbf{H}^1(\mathcal{O})} & \leq Ch|\tilde{f}(\lambda w_1 - w_1^*) - \gamma_0^+(\lambda w_1 - w_1^*)|_{2, \mathcal{O}} \notag \\
& \leq Ch\|\lambda w_1 - w_1^*\|_{H^2(\Omega)} \label{uhEst3}
\end{align}

By the continuity of the right inverse of $\gamma_0 : \mathbf{H}^1(\mathcal{O}) \to \mathbf{H}^{1/2}(\partial \mathcal{O})$ and the estimate
(\ref{w1Est1}) or (\ref{w1Est2}), we also have
\begin{align}
\|\gamma_0^+(\lambda w_{1} - w_1^*) - \gamma_0^+(\lambda w_{1h} - w_1^*)\|_{\mathbf{H}^1(\mathcal{O})} &\leq C_\lambda \|w_1 -
w_{1h}\|_{H^1(\Omega)} \notag\\
&\leq C_\lambda h \big\|[w_1^*, w_2^*, u^*]\big\|_{\mathbf{H}_\rho}. \label{uhEst4}
\end{align}

From (\ref{map1}) and (\ref{muEstimate}) as well as either (\ref{w1Est1}) or (\ref{w1Est2}) we also have
\begin{align}
\|\tilde{f}_{0h}(\lambda w_{1} - w_1^*) - \tilde{f}_{0h}(\lambda w_{1h} - w_1^*)\|_{\mathbf{H}^1(\mathcal{O})} &\leq C_\lambda \|w_1 -
w_{1h}\|_{H^2(\Omega)} \notag\\
&\leq C_\lambda h \big\|[w_1^*, w_2^*, u^*]\big\|_{\mathbf{H}_\rho}. \label{uhEst5}
\end{align}

Applying (\ref{uhEst2})-(\ref{uhEst5}) to the right hand side of (\ref{uhEst1}) (as well as considering the continuous dependence of the data
inherent in Theorem \ref{well}$(i)$,) we then have
\begin{equation}
\|u - u_h\|_{\mathbf{H}^1(\mathcal{O})} \leq C_\lambda h \big\|[w_1^*, w_2^*, u^*]\big\|_{\mathbf{H}_\rho},
\end{equation}
which proves Theorem \ref{thm:estimates} $(iii)$.

Finally, for the error in the fluid pressure term: by (\ref{up}) and (\ref{ph}) we have
\begin{align}
\|p - p_h\|_{L^2(\mathcal{O})} & \leq \|\tilde{\pi}(\lambda w_1 - w_1^*) - \tilde{\pi}_h(\lambda w_{1h} - w_1^*)\|_{L^2(\mathcal{O})}
 + \|\tilde{q}(u^*) - \tilde{q}_h(u^*)\|_{L^2(\mathcal{O})} + C_\lambda h \big\|[w_1^*, w_2^*, u^*]\|_{\mathbf{H}_\rho}
\end{align}
after using (\ref{cEst4}).  Proceeding just as in the proof of Theorem \ref{thm:estimates}$(iii)$ above results in the asserted estimate
(\ref{ivpressure}). This completes the proof of Theorem \ref{thm:estimates}.   \ \ \ $\square $

\subsection{Matlab Implementation of the Finite Element Method} We include here a brief description of the numerical scheme followed by some
numerical results from a test problem.  The finite element method is a numerical implementation of the Ritz-Galerkin method over a specific set
of basis functions defined on a mesh of the domain.  In this case the fluid domain is divided into tetrahedra and the plate domain into
triangles. Basis functions are associated to points in the mesh and the system is solved in this finite dimensional setting via a matrix/vector
equation, see e.g. \cite{Axelsson}.

First consider the plate system (\ref{BB1}).  This weak formulation takes place over $H^2_0(\Omega)$ and thus the most natural choice for the
discretization is a set of $H^2$-conforming elements, see \cite{Solin}; for a MATLAB implementation see \cite{Dominguez}. We use the quintic
Argyris basis functions because they are the lowest order $H^2_0$-conforming elements available and they ensure wellposedness of the discrete
formulation of (\ref{BB1}) because the inf-sup condition on the bilinear form $b$ can still be satisfied. The Argyris basis functions have 21
degrees of freedom(DOF) for each triangle in the mesh, namely Lagrange DOF for function values at each vertex, Hermite DOF for
$\partial/\partial x$ and $\partial/\partial y$ at each vertex, Argyris DOF for $\partial^2/\partial x^2$, $\partial^2/\partial x\partial y$,
and $\partial^2/\partial y^2$ at each vertex and one DOF at each edge midpoint for the normal derivative.

On the fluid domain a Stokes system is solved via a mixed formulation.  Again one must be careful as to the choice of elements used to guarantee
that the discrete problem remains well-posed. In this case we use the popular Taylor-Hood $(\mathbb{P}^2/\mathbb{P}^1)$ elements, see e.g.
 \cite{Ern}.

To derive the discretized version of (\ref{BB1}) let $X_h$ denote the space of Argyris basis functions of dimension $n$ with basis
$\{\phi_i\}_{i=1}^n$.  Let $w_h = \sum_{i=1}^n \alpha_i \phi_i$.  Then for each $\phi_j \in X_h$ and $r\in \mathbb{R}$ we have
\begin{align*}
\sum_{i=1}^n \alpha_i a_\lambda(\phi_i, \phi_j) + b(\phi_j,\tilde{c}) &= \mathbb{F}(\phi_j);\\
\sum_{i=1}^n \alpha_i b(\phi_i,r) &=0,
\end{align*}
with $a_\lambda$, $b$ and $\mathbb{F}$ defined in (\ref{BB1}) - (\ref{abF}).  By writing this for each $\phi_j \in X_h$  we build a linear
system of the form
\begin{equation} \begin{bmatrix} A & B \\B^T & 0\end{bmatrix}
\begin{bmatrix} \alpha
\\ \tilde{c} \end{bmatrix} = \begin{bmatrix} \mathbf{F}\\0\end{bmatrix}. \label{MV1}
\end{equation}
Note importantly, that in practice $\tilde{f}$ and $\tilde{\mu}$ are not known exactly and so in reality $(\ref{MV1})$ is actually created using
a subroutine that numerically approximates the instances of $\tilde{f}$ and $\tilde{\mu}$ that occur in $a_\lambda$ using a discrete mixed
variational formulation on the fluid domain. The discretized versions of (\ref{sol1}) and (\ref{sol2}) take the same mixed form as (\ref{MV1}),
but with $\tilde{\mathbf{a}}(\mu,\varphi)$ and $\tilde{\mathbf{b}}(\mu,q)$ defined as in (\ref{atilde}) - (\ref{btilde}).

Here we consider a test problem for the fluid-structure problem of interest in the $\rho = 0$ case.  The fluid domain $\mathcal{O}$ is given by
$(x,y,z) \in (0,1) \times (0,1) \times (-1,0)$.  The plate $\Omega$ is the top boundary of the fluid domain, lying in the $xy$ plane, namely
$(x,y)\in (0,1)\times(0,1)$.

Then for any $\lambda >0$ the functions
\begin{align*}
w_1 &= -x^4(x-1)^4 (2x-1) y^4 (y-1)^4 ;\\
w_2 &= -\Delta w_1 = 12x^2(x-1)^2(2x-1)(6x^2-6x+1) y^4 (y-1)^4 \\
& \hspace{.8in}+ x^4(x-1)^4 (2x-1) 4y^2 (y-1)^2 (14y^2-14y+3);\\
u^1 &= \big[2x^3(x-1)^3(9x^2-9x+2)y^4(y-1)^4 + (4/5)x^5(x-1)^5y^2(y-1)^2(14y^2-14y+3)\big]*\\
&\hspace{.35in}\big[-30z^4-60z^3-30z^2\big] ;\\
u^2 &= 0 ;\\
u^3 &= -\big[12x^2(x-1)^2(2x-1)(6x^2-6x+1)y^4(y-1)^4 + 4x^4(x-1)^4(2x-1)y^2(y-1)^2(14y^2-14y+3)\big]*\\
&\hspace{.35in}\big[-6z^5 - 15z^4 - 10z^3 -1\big];\\
p &=0,
\end{align*}
solve (\ref{w2}) - (\ref{end}) for data defined by
\begin{align*}
w_1^* &\equiv \lambda w_1 - w_2; \\
w_2^* &\equiv \lambda w_2 + \Delta^2 w_1;\\
u^* &\equiv \lambda u - \Delta u.
\end{align*}

Notice that $u = [u^1,u^2,u^3]$ is divergence free, $u = [0,0,0]$ on $S$, $\Delta u \cdot \nu = [0,0,0]$ on $S$, and $u=[0,0,w_2]$ on $\Omega$.
Moreover $p:=G_{\rho,1}(w_1) + G_{\rho,2}(u) = 0$ because $w_1$ and $u$ are chosen such that $\Delta^2 w_1 = - \Delta u \cdot \nu$ on $\Omega$,
and $\Delta u \cdot \nu = 0$ on $S$ which causes the two terms to cancel. Finally we have $[w_1^*, w_2^*, u^*] \in \mathbf{H}_\rho$ (in the
$\rho = 0$ case $w_2^* \in L^2(\Omega)/\mathbb{R}$ only).  The error in the numerical solution for this test problem is summarized in the table
below.
\begin{table}[H]
\centering
\begin{tabular}{c c c c c c}
\hline
No. of elements & Characteristic Length & $|w_1 - w_{1h}|_{H^2}$ & $|w_1 - w_{1h}|_{H^1}$ & $\|w_1 - w_{1h}\|_{L^2}$\\
\hline
4  &  1 & $7.132 \times 10^{-5}$ & $4.993 \times 10^{-6}$ & $ 3.935\times 10^{-7}$ \\
16 & .5 & $9.823\times 10^{-6}$ & $3.450\times 10^{-7}$ & $1.509\times 10^{-8}$ \\
64 & .25 & $1.249\times 10^{-6}$ & $2.761\times 10^{-8}$ & $1.598\times 10^{-9} $    \\
256 & .125 & $8.343\times 10^{-8}$ & $1.253\times 10^{-9}$ & $1.066\times 10^{-10} $   \\
1024 & .0625 & $5.124\times 10^{-9}$ & $6.771\times 10^{-11}$ & $6.285\times 10^{-12} $   \\
\hline
\end{tabular}
\caption{Errors of structure FEM approximations.}
\end{table}
Since the mesh is refined by a factor of 2 at each step, we compute $\log\Big(\frac{\text{Error}_i}{\text{Error}_{i+1}}\Big)/\log(2)$.  In the
limit this ratio should approach the exponent of convergence, i.e. $\mathcal{O}(h^k)$.  Now for smooth data (as we have here) the best possible
convergence rate one could attain is $k=4$ for the $H^2$ norm of $w_1$ which the numerical scheme does appear to attain (see Table
\ref{tab:w1}). However, the $H^1$ and $L^2$ errors do not appear to improve to $k=5$ and $k=6$ respectively; this is possibly due to the
(unavoidable) approximation of $\tilde{f}$ and $\tilde{\mu}$ described above.
\begin{table}[H]
\centering
\begin{tabular}{c c c c}
\hline
& $H^2$ & $H^1$ & $L^2$\\
\hline %
Mesh 1 / Mesh 2 & 2.86 & 3.86   & 4.71 \\
Mesh 2 / Mesh 3 & 2.98 & 3.64 & 3.24 \\
Mesh 3 / Mesh 4 & 3.90 & 4.46 & 3.91 \\
Mesh 4 / Mesh 5 & 4.03 & 4.21 & 4.08 \\
\hline
\end{tabular}
\caption{Computed index $k$ in $\mathcal{O}(h^k)$ for structure FEM approximations.} \label{tab:w1}
\end{table}

\begin{figure}[H]
\centering
  \includegraphics[width=.6\linewidth]{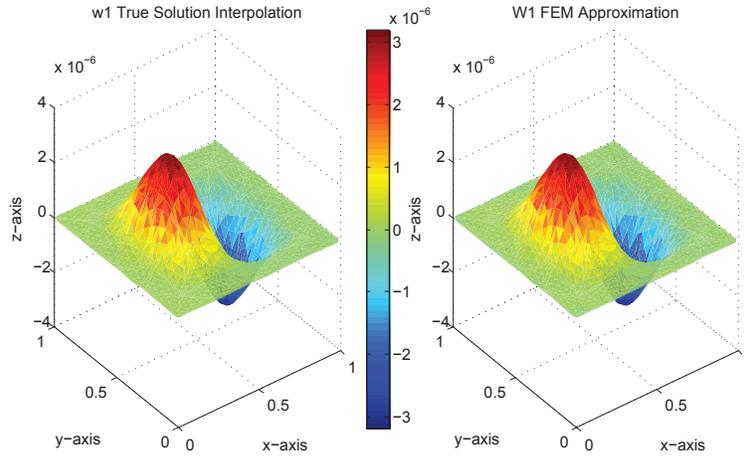}
  \caption{FEM approximation for $w_1$.}
  \label{fig:w1}
\end{figure}

In Figure \ref{fig:w1} we see that already at the 3rd level of mesh the FEM approximation is indistinguishable from the true solution.

Similarly, for the fluid approximation we have in Table \ref{tab:fluid} the errors in the fluid variables for each mesh refinement.
\begin{table}[H]
\centering
\begin{tabular}{c c c c c c}
\hline
No. of elements & Characteristic Length & $\|u - u_{h}\|_{\mathbf{L}^2}$ & $|u - u_{h}|_{\mathbf{H}^1}$ & $\|p - p_{h}\|_{L^2}$\\
\hline %
24 & 1 & $5.26\times 10^{-4}$ & $9.53\times 10^{-3}$ & $1.40\times 10^{-4}$\\
192 & .5 & $1.42\times 10^{-5}$ & $3.64\times 10^{-4}$ & $5.85\times 10^{-5}$\\
1536 & .25 & $3.56\times 10^{-6}$ & $1.75\times 10^{-4}$ & $2.25\times 10^{-5}$\\
12288 & .125 & $4.98\times 10^{-7}$ & $5.20\times 10^{-5}$ & $3.56\times 10^{-6}$\\
98304  & .0625 & $6.33\times 10^{-8}$ & $1.37\times 10^{-5}$ & $6.67\times 10^{-7}$\\
\hline
\end{tabular}
\caption{Errors of fluid FEM approximations.} \label{tab:fluid}
\end{table}
The log error ratios approach what is expected for a $\mathbb{P}^2/ \mathbb{P}^1$ implementation, namely $k = 3, 2,$ and $2$ respectively as
shown in Table \ref{tab:FluidRate}.
\begin{table}[H]
\centering
\begin{tabular}{c c c c}
\hline
& $\mathbf{L}^2$(fluid) & $\mathbf{H}^1$(fluid) & $L^2$(pressure)\\
\hline %
Mesh 1 / Mesh 2 & 1.89 & 1.38 & 1.26 \\
Mesh 2 / Mesh 3 & 1.99 & 1.05 & 1.38 \\
Mesh 3 / Mesh 4 & 2.84 & 1.75 & 2.66 \\
Mesh 4 / Mesh 5 & 2.97 & 1.93 & 2.41 \\
\hline
\end{tabular}
\caption{Computed index $k$ in $\mathcal{O}(h^k)$ for fluid FEM approximations.}
\label{tab:FluidRate}
\end{table}

\begin{figure}[H]
\centering
  \includegraphics[width=.6\linewidth]{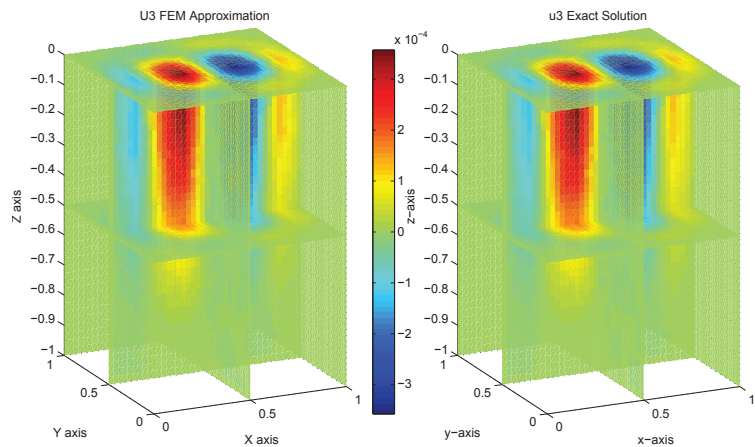}
  \caption{FEM approximation for $u^3$.}
  \label{fig:u3}
\end{figure}

\begin{figure}[H]
\centering
  \includegraphics[width=.6\linewidth]{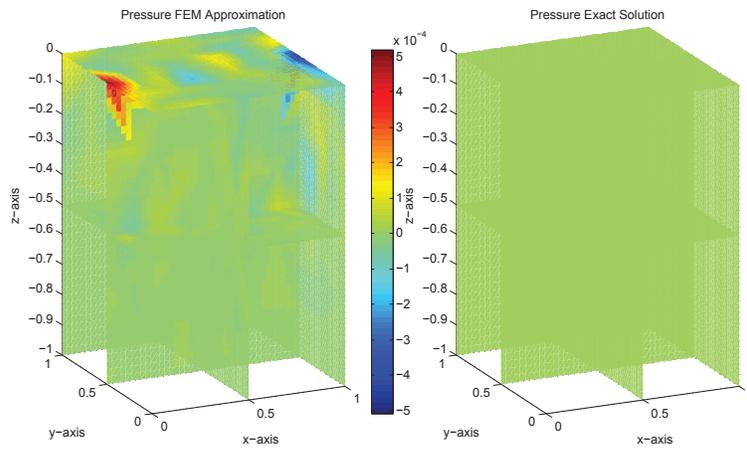}
  \caption{FEM approximation for $p$.}
  \label{fig:p}
\end{figure}

In Figure \ref{fig:u3} we see a slice of the third component of the 3-D fluid velocity $u^3$ displaying the test problem's non-trivial boundary
interaction with the plate. Figure \ref{fig:p} shows the pressure is converging to the solution as well.

\begin{figure}[H]
\begin{center}
\includegraphics[width=7cm]{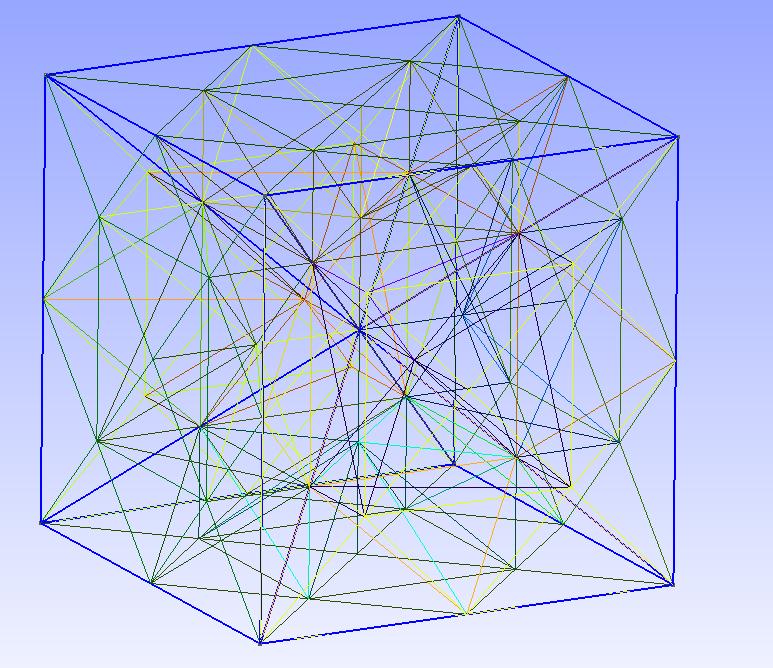}
\end{center}
\caption{3D fluid mesh created with GMSH.} \label{GMSH}
\end{figure}

\end{document}